\documentclass[11pt]{elsart}
\usepackage{graphicx}
\usepackage{subfigure}
\usepackage{multirow}
\usepackage{url}
\usepackage{amssymb}
\usepackage{amsmath}
\usepackage{amsfonts}
\usepackage{wrapfig}
\usepackage{afterpage}
\usepackage{float}
\usepackage{algorithm}
\usepackage[noend]{algpseudocode}
\usepackage{epsfig}
\usepackage{epstopdf}

\usepackage{url}

\graphicspath{{fg/}}
\newcommand{\real}{\mathbb{R}}
\newcommand{\beps}{\mbox{\boldmath$\epsilon$}}
\newcommand{\brho}{\mbox{\boldmath$\rho$}}
\newcommand{\bsig}{\mbox{\boldmath$\sigma$}}
\newcommand{\btau}{\mbox{\boldmath$\psi$}}
\def\vsig{\varsigma}

\newcommand{\bveps}{\mbox{\boldmath$\varepsilon$}}
\newcommand{\bxi}{\mbox{\boldmath$\xi$}}

\newcommand{\bvsig}{\mbox{\boldmath$\varsigma$}}

\newcommand{\bzeta}{\bvsig} 

\def\bw{{\bf w}}
\def\uu{u}
\def\WW{W}
\def\ww{w}
\def\calW{{\cal W}}

\def\la{\langle}
\def\ra{\rangle}

\def\bR{{\bf R }}
\def\calF{{\cal F}}
\def\FF{{F}}
\def\alp{\alpha}
\def\lam{\lambda}

\def\vv{v}

\def\bz{\bf z}

\def\bb{{\bf b}}
\def\half{\frac{1}{2}}
\def\t{\mathbf{t}}

\def\u{\mathbf{u}}

\def\Diag{\mbox{Diag}}

\def\bG{\mathbf{G}}
\def\bF{\mathbf{F}}
\def\calS{{\cal S}}

\def\calP{{\cal P}}
\def\barbsig{\bar{\bsig}}
\def\bartau{\bar{\tau}}
\def\barbrho{\bar{\brho}}
\def\barbu{\bar{\u}}
\def\eb{\begin{equation}}
\def\ee{\end{equation}}
\newtheorem{exam}{Example}

\def\PP{P}

\def\bx{\mathbf{x}}
\def\by{\mathbf{y}}

\def\eps{\epsilon}
\def\sig{\sigma}

\def\calX{{\cal X}}

\def\aa{v}
\def\cc{w}
\def\Lam{\Lambda}
\def\calZ{{\cal Z}}
\def\calU{{\cal U}}

\def\bb{{\bf b}}
\def\half{\frac{1}{2}}
\def\t{\mathbf{t}}
\def\bt{\mathbf{t}}

\def\u{\mathbf{u}}
\def\bu{\mathbf{u}}

\def\bv{\mathbf{v}}

\def\bc{\mathbf{c}}
\def\bD{\mathbf{D}}
\def\bK{\mathbf{K}}

\def\Diag{\mbox{Diag}}

\def\bC{\mathbf{C}}
\def\bQ{\mathbf{Q}}

\def\ot{ \circ}

\def\bG{\mathbf{G}}
\def\bF{\mathbf{F}}
\def\bff{\mathbf{f}}
\def\calS{{\cal S}}
\def\calE{{\cal E}}

\def\calP{{\cal P}}

\def\bN{\mathbf{N}}
\def\dO{{{\rm d}\Omega}}
\def\dG{{{\rm d}\Gamma}}
\def\barbsig{\bar{\bsig}}
\def\barbvsig{\bar{\bvsig}}

\def\barvsig{\bar{\tau}}

\def\barbx{\bar{\bx}}

\def\barbrho{\bar{\brho}}
\def\barbzeta{\bar{\bzeta}}
\def\barbu{\bar{\u}}
\def\eb{\begin{equation}}
\def\ee{\end{equation}}
  \newtheorem{definition}{Definition}
 \newtheorem{remark}{Remark}

\begin{document}
\begin{frontmatter}
\title{On Topology Optimization and \\ Canonical Duality  Method}
\author {David Yang Gao\corauthref{cor1}}
\address {School  of Science and Technology,  Federation University Australia}
  \corauth[cor1]{Corresponding author: d.gao@federation.edu.au}


\begin{abstract}
Topology optimization for general materials  is correctly formulated as
  a   bi-level knapsack problem,
 which is considered to be NP-hard in global optimization and  computer science.
By using   canonical duality theory   (CDT)  developed  by the author, the linear  knapsack  problem can be solved analytically to obtain global optimal solution at each design iteration.  Both uniqueness, existence,  and NP-hardness  are discussed.
  The novel  CDT method for general topology optimization  is
  refined and  tested by both 2-D and 3-D benchmark problems.
Numerical results show that  without using  filter and any  other artificial technique, the   CDT method can produce  exactly 0-1 optimal density distribution with almost no  checkerboard pattern.
   Its performance and novelty are compared with   the popular  SIMP and  BESO approaches.
    Additionally, some  mathematical and conceptual mistakes  in literature are explicitly  pointed out.
   A brief review on the canonical duality theory for solving a unified problem in multi-scale nonconvex/discrete  systems is given in Appendix.
\end{abstract}

\end{frontmatter}

\section{Introduction}
Topology optimization is a mathematical tool  that optimizes material layout within a  prescribed design domain in order to  obtain the best structural
performance  under    a given set of loads and  geometric/physical constraints.
 Due to its broad applications, this  tool  has   been studied extensively by both engineers and mathematicians for more than 40 years   (see the   comparative reviews \cite{bb_Rozvany2009,sig-mau}).
 Generally speaking, a topology optimization problem  involves both continuous state variable (such as the deformation field ${\bf u(x)}$) and the  density distribution $\rho({\bf x})$ that can take either the value 0 (void) or 1 (solid) at any point in the design domain. Thus,  numerical discretization methods  (say FEM) for solving  topology optimization problems lead to a so-called  mixed integer nonlinear programming (MINLP)  problem, which appears not only in   computational engineering, but also in operations research, decision and management sciences,  industrial and systems engineering \cite{gao-opl16}.

As one of the most challenging problems in global optimization, the MINLP has been seriously studied by mathematicians and computational scientists for several decades,
many  methods and algorithms have been proposed (see~\cite{Lee-Ley}). These methods can be categorized into two main groups \cite{g-r-l}: deterministic   and stochastic methods. The stochastic methods are based on an element of random choice. Because of this, one has to sacrifice the possibility of an absolute guarantee of success within a finite amount of computation. The deterministic methods, such as  cutting plane, branch and bound, semi-definite programming (SDP), can find global optimal solutions, but not in polynomial time. Therefore,   the MINLP is known  to be NP-hard (non-deterministic polynomial-time hard). Indeed, even the most simple quadratic integer 0-1 programming
\[
\min_{{\bx } \in {\mathbb R}^n}  \left\{ \frac{1}{2} \bx^T \bQ \bx - \bx^T \bff | \;\;    \bx  \in \{ 0,1\}^n \right\}
\]
is considered to be  NP-hard. This integer minimization problem has
$2^n$  local solutions. Due to the lack of global optimality criteria, traditional direct approaches  can only handle very small size problems.
Therefore,  global optimization  problems with 200 variables
are referred to as ``medium scale",  problems with 1,000 variables as ``large scale", and the so-called
``extra-large scale" is only around 4,000 variables \cite{bura}.
It was proved by Pardalos and Vavasis  \cite{p-v}  that instead of the integer constraint, even  the continuous
 quadratic minimization with box constraints $ 0 \le \bx \le 1 $  is NP-hard as long as the
  matrix $\bQ $ has one negative eigenvalue. However,   it was discovered by the author that
these   so-called NP-hard problems can be solved easily by {\em canonical duality theory}  as long as the  global optimal solution is unique \cite{gao-jimo07,gao-cace,gao-opl16}.

By the fact that the topology optimization has to handle huge-scale   MINLP problems with   millions of variables,
various relaxation  approaches  and  techniques  have been developed by engineers, such as  the homogenization~\cite{bb_Bendsoe88}, density-based method~\cite{bb_Bendsoe89},  phase field approach,  topological derivatives~\cite{bb_Sokolowski99,bb_Suresh2010A}, the
level set methods ~\cite{bb_Allaire2002A,bb_Wang2003A}, as well as the well-known   SIMP  (Simplified Isotropic Material with Penalization)   \cite{bb_Zhou1991} and
 evolutionary methods (ESO and BESO) ~\cite{bb_Huang20102,bb_Mattheck1990A,xie93}.
Most of these  engineering approaches generally relax the MINLP   as  a continuous parameter optimization problem, and  then solve it based on
the traditional   Newton-type (gradient-based) methods.
There exists several fundamental issues on these approaches. First,
 the relaxation from discrete  to continuous optimization must be  mathematically correct. Otherwise, the
 numerical results produced by these methods can't convergent to  mechanically sound structural topology.
Second,  the Newton-type algorithms can be used only for convex minimization. For nonconvex problems,
  numerical results obtained by these algorithms depend    sensitively on  initial data and numerical precisions adopted.
  It was discovered  by Gao  and his co-workers
   that the global optimal solutions are usually nonsmooth not only  for coupled optimal design problems (see Chapter 7, \cite{gao-dual00}), but also for general
   nonconvex variational problems \cite{gao-ogden-qjmam}. These nonsmooth  solutions can't be captured by any Newton-type algorithms.
By the fact that the  SIMP  is not a mathematically correct penalty method,  this most popular engineering approach  can never produce exact integer solution for any given real-world problem.
The  existence of gray scale elements and appearance of checkerboards patterns are the SIMP's  two major  intrinsic problems
\cite{bb_Diaz1995,bb_Sigmund1998}.
 Although the  commercial code by  the  BESO can produce integer solutions, it was discovered recently \cite{gao-ali-to} that  this popular method is actually a direct approach for
 solving a knapsack-type problem and it is not
 a  polynomial-time algorithm.   This the reason why the  BESO is computationally expensive and  can be used only for small sized problems.

Duality approaches for topology optimization have been studied via the traditional Lagrangian duality theory  \cite{bb_Beckers1999,bb_Beckers2000,bb_Jog2001,bb_Jog2002,bb_Jog2010A,bb_Stolpe2003}.
However,  the Lagrange multiplier method can be used mainly for convex problems with equality constraints~\cite{l-g-opl}.
 For nonconvex problems, the Lagrangian $L(\bx, \by)$ is usually not a saddle function.
By the fact that
\[
\min_{\bx} \max_{\by} L(\bx, \by) \ge \max_{\by} \min_{\bx} L(\bx, \by) ,
\]
the Lagrangian duality theory produces a so-called duality gap at each iteration. In order to reduce this duality gap, much effort has been made  by mathematicians during  the past 30 years~\cite{gao-mot,Manyem}. For inequality constraints, both the Lagrange multiplier and the constraint must satisfy the KKT conditions. The associated complementarity problem is very difficult even for linearly  constrained problems in continuous space~\cite{isac}.
Although the augmented Lagrange multiplier  method can be used for solving both equality and  inequality constrained problems,   the   constraints must be
   linear     since any   simple  nonlinear constraint could lead to a nonconvex minimization problem \cite{l-g-opl}.
Unfortunately, all these mathematical difficulties were not correctly addressed   in the topology optimization literature (see \cite{bb_Jog2001,bb_Jog2002,bb_Jog2010A}).

Canonical duality theory  (CDT) is a  precise  methodological theory, which can be used not only for modeling complex systems within a unified framework \cite{gao-aip}, but also for solving a large class of challenging problems in nonconvex analysis and global optimization \cite{gao-dual00,g-18,g-l-r-17}.  Application of this theory to general topology optimization was given recently \cite{g-to,g-18}. It was  discovered  that   the  0-1 integer programming in  topology optimization for linear elasticity
 is actually equivalent to the well-known Knapsack problem, which can be solved analytically by the CDT \cite{gao-jimo07}.
The main goal of this paper is to present a detailed  study on the canonical duality approach  for solving general topology optimization
problems  with applications to 2-D and 3-D linear elastic structures.
In the next section, the general topology optimization problem and its challenges are  discussed.
A mathematically correct topology optimization problem is  formulated as a coupled bilevel knapsack problem.
The conceptual mistakes in   topology optimization and mathematical difficulties in SIMP method are discussed.
 Section~\ref{sec-dual} shows that the knapsack problem can be solved analytically to obtain global optimal solution at each iteration.
 A canonical dual  algorithm  for computing the globally optimal dual solution is explained in Section~\ref{sec-num}.
 Some fundamental issues on  challenges in topology optimization and  NP-hardness in computational complexity are addressed in Section \ref{sec4}.
 Numerical examples are shown  in Section~\ref{sec-example}. Conclusion remarks are given in Section~\ref{sec-con}.
 A  brief review of the canonical duality theory is provided  in Appendix.

\section{On Mathematical Models   and Challenges}\label{sec-basics}
Let us consider an elastically deformable body that in an undeformed configuration occupies an open  domain $\Omega\subset \real^d \; (d=2,3)$ with boundary $\Gamma = \partial \Omega $. We assume that the body is subjected to a body force $\bb$ (per unit mass) in the reference domain $\Omega$  and a given surface traction $\bt(\bx)$ of dead-load type  on the boundary ${\Gamma_t} \subset \partial \Omega $, while the body is fixed on the remaining boundary ${\Gamma_u} = \partial  \Omega\setminus {\Gamma_t}$.
The total potential energy of this deformed body has the following standard form:
\begin{equation}
 \Pi(\bu, \rho) = \int_{\Omega}
W(\nabla \bu) \rho {\rm d} \Omega - \int_{\Omega} \bu \cdot \bb   \rho
{\rm d} \Omega - \int_{\Gamma_t} \bu \cdot \bt {\rm d} \Gamma
  , \label{pprobm}
\end{equation}
where the displacement  $\bu : \Omega \rightarrow \real^d$ is a continuous field variable,  the mass density  $\rho:\Omega\rightarrow
 \{ 0, 1\}$ is a discrete  design variable;
the stored  energy density  $W(\bF)$ is an {\em objective function} (see Appendix)   of the deformation gradient $\bF=\nabla \bu$.
It should be emphasized that   the design variable $\rho$  is in a discrete space subjected to  a
so-called {\em  knapsack condition}  $\int_\Omega \rho(\bx) \dO \le V_c $ ($V_c > 0$ is a given  desired  volume),
   classical  variational method can't be used  to obtain the criticality condition for $\rho$.
Therefore, analytical methods for  studying general  topology optimization problems are fundamentally difficult.

By using finite element method, the domain $\Omega$ is divided into $n$
 elements $\{\Omega_e\}$ and in each $\Omega_e$, the unknown fields can be numerically discretized as
\begin{equation}\label{eq-urho}
\bu(\bx) = \bN_e(\bx) \bu_e , \;\; \rho(\bx) = \rho_e \in \{0,  1 \} \;\; \forall \bx \in \Omega_e,
\end{equation}
where $\bN_e$ is an interpolation matrix, $\bu_e$ is a nodal displacement vector, the binary design variable $\rho_e \in \{ 0,1\}$
 is used for determining whether the element $\Omega_e $ is a void ($\rho_e = 0$) or a solid ($\rho_e = 1$).
Let  $\calU_a \subset \real^m$   be a kinetically  admissible  nodal  displacement space,
 $\vv_e > 0$ be the volume of the $e$-th element $\Omega_e$,  and
\eb
 \calZ_a =
 \left \{ \brho = \{ \rho_e \} \in   \{ 0, 1\}^n |
  \;\;  \brho^T\bv = \sum_{e=1}^n \rho_e \vv_e \le V_c \right\}.
\ee
Thus, by substituting~\eqref{eq-urho} into   \eqref{pprobm}, the total potential energy functional
can be numerically reformulated as  a real-valued function $\Pi_h:\calU_a \times \calZ_a \rightarrow \real$:
\begin{equation}\label{eq-ph}
   \Pi_h(\bu, \brho) = C(\brho,\bu)  - \bu^T \bff  ,
\end{equation}
where  $C(\brho, \bu) = \brho^T \bc(\bu)$,
\eb
\bc(\bu) =  \left\{
\int_{\Omega_e} [ W(\nabla \bN_e(\bx) \bu_e) -   \bb_e^T \bN_e(\bx)\bu_e ] \dO   \right\}
 \in \real^n,
\label{eq-cu}
 \ee
\eb
\bff  =   \left\{ \int_{\Gamma^e_t} \bN_e (\bx)^T \t_e (\bx) \dG \right\}  \; \in \real^m.
\ee

 By the facts that $\brho \in \calZ_a$ is the main design variable,
 the displacement $\bu$ depends on each given domain $\Omega$, and
$\brho^p = \brho \;\; \forall p\in\real, \;\; \brho \in \calZ_a$,
  the   topology optimization for  general  nonlinearly  deformed structures
 should  be formulated as a  so-called bi-level   mixed integer nonlinear programming \cite{g-18}:
\begin{eqnarray}
 (\calP_{bk}):\;\; &  \;\;\;\;\;\; &  \min_{\brho\in \calZ_a}
\{  \Phi_p (\brho,\barbu ) = \bff^T \barbu   -\bc(\barbu)^T \brho^p  \}  \label{eq-upp}  \\
 & & \mbox{s.t.}  \;\; \barbu(\brho)  =\arg \min_{\u \in \calU_a } \Pi_h(\bu, \brho). \label{eq-llopt}
\end{eqnarray}
In this formulation, $\Phi_p (\brho,\bu)   $ represents the upper-level cost function and the total potential energy $\Pi_h(\bu, \brho)$ represents the lower-level cost function.
By the fact that for each given $\barbu$, the upper-level optimization is actually a typical  knapsack problem,
$(\calP_{bk})$ is essentially a coupled bilevel knapsack problem (BKP).
Bilevel programming is known to be strongly NP-hard \cite{h-j-s}, and it has
been proven that merely evaluating a solution for optimality
is also a NP-hard task \cite{v-s-j}. Even in the simplest case of
linear bilevel programs, where the lower level problem has a
unique optimal solution for all the parameters, it is not likely
to find a polynomial algorithm that is capable of solving the
linear bilevel program to global optimality \cite{s-m-d}. The proof for the
non-existence of a polynomial time algorithm for linear bilevel
problems can be found in \cite{deng}.
For large deformation problems, the total potential energy $\Pi_h$  is usually a nonconvex function of $\bu$.
  Therefore, this bi-level optimization should  be the most challenging problem so far  in global optimization and computational mechanics.

In reality,  the topology optimization is a design process, an alternative iteration method can be naturally
used to solve the bilevel optimization problem, i.e.
  \begin{verse}
  (i) For a given $\brho_{k-1} \in \calZ_a$, to solve the lower-level problem (\ref{eq-llopt})  first for
  \eb
  \bu_k =  \arg \min \{ \Pi_h( \bu, \brho_{k-1})   \;\; | \;\; \bu\in \calU_a \} . \label{eq-lowk}
  \ee

  (ii) Then, for the fixed  $\bu_k \in \calU_a$, to solve the upper-level integer minimization problem (\ref{eq-upp}) for   $\brho_k\in \calZ_a$ such that
  \eb
\brho_k = \arg   \min \{ \Phi_p  (\brho, \bu_k )  | \; \;\; \bz\in \calZ_a \} \label{eq-uppk}.
  \ee
  \end{verse}
The canonical duality and finite element method for solving \eqref{eq-lowk} has been studied extensively for computational mechanics and global optimization \cite{gao-jem96,g-l-r-17}.
This paper will focus only on the   upper-level problem \eqref{eq-uppk}.

 Let   $\bw = \bc(\bu) \in \real^n  $ and $p=1$.
The  upper-level optimization  \eqref{eq-uppk}   is equivalent to a   standared linear knapsack problem  ($(\calP_{kp}) $ for short) \cite{g-to}:
\eb\label{eq-low}
 (\calP_{kp}): \;\;  \min\left\{  \PP_u(\brho) = -   \bw^T   \brho  \;\; 
  |  \;\; \bv^T \brho \le V_c, \;\; \brho \in  \{0, 1\}^n  \right\}.
\ee
This well-known  problem in  decision science makes a perfect sense in topology optimization, i.e.
among all elements $\{\Omega_e\}$, one should keep only those who stored more deformation  energy $ \{ w_e \} = \{ c_e(\bu)\}$.
Due to the integer constraint, even this linear 0-1 programming is listed as one of Karp's 21 NP-complete problems~\cite{karp}. However, this challenging problem can be solved analytically by using the canonical duality theory.

For   linear elastic structures without body force,   the total potential energy is simply a quadratic function:
\eb
 \Pi_h(\bu, \brho) = \half \bu^T \bK(\brho) \bu - \bu^T \bff   ,  \label{eq-lep}
\ee
where $ \bK(\brho) = \left\{ \rho_e \bK_e \right\} \in \real^{n\times n} $ is the overall stiffness matrix, obtained by assembling the sub-matrix $\rho_e \bK_e$ for each element $\Omega_e$.
For a given $\brho \in \calZ_a$,  the global optimal solution for the lower-level minimization problem
(\ref{eq-llopt})   is  governed  by
 $ \bK(\brho) \barbu =  \bff .$
In this case,  $\bc(\bu) = \half \left\{  \bu^T_e  \bK_e  \bu_e \right\} \in \real^n_+
= \{ \bc \in \real^n | \;\; \bc \ge {\bf 0} \in \real^n \} $.
Then the topology optimization  problem  for linear elastic structures   can be   written in the following form
\begin{eqnarray}
 (\calP_{le}):\;\; &     & \min_{\brho\in \calZ_a}
\{   \bff^T \barbu(\brho) - \brho^T \bc( \barbu(\brho))    \}   \label{eq-le} \\
&  \mbox{s.t.} &   \barbu(\brho)  =\arg \min_{\u \in \calU_a }  \left\{   \half \bu^T \bK(\brho) \bu - \bu^T \bff  \right\}\label{eq-lel}
\end{eqnarray}
 This is a typical bilevel knapsack-quadratic  optimization.
By the fact that the lower-level problem has a unique solution governed by $\bK(\brho) \bu = \bff$,
the single-level reduction for solving  this problem leads to
\eb
 (\calP_{sl}):\;\;   \min_{\brho\in \calZ_a}
\left\{ \bff^T \bu - \half \bu^T \bK( \brho)  \bu      \;  | \;\;   \bK(\brho) \bu = \bff , \;\; \bu \in \calU_a
  \right \} .
\ee
\begin{remark}[On Minimum Compliance Problem and SIMP Method]
{\em
Instead of $(\calP_{le})$ or $(\calP_{sl})$,   the topology optimization problem in literature is usually formulated as
  the  so-called minimum compliance problem \cite{bb_Bendsoe89,sig-mau}:
\eb
 (\calP_c): \;\;\; \min
\left\{ \half \bu^T \bff  \;\;  |  \;\;   \bK(\brho) \bu =  \bff ,     \;\;  \brho \in \calZ_a  \right\},  \label{eq-mcp}
\ee
where  the linear cost $  \half \bu^T  \bff  $   is called  the ``mean compliance''  in topology optimization \cite{sig-mau}.
If the state variable  is replaced by $\bu = \bK(\brho)^{-1} \bff$, then $(\calP_c)$ can be written as
\eb
 (\calP^c_{sl}):\;\;   \min_{\brho\in \calZ_a}
\left\{  \half \bff^T [ \bK( \brho)  ]^{-1} \bff    \;  | \;\;   \bK(\brho) \mbox{ \rm is invertible on } \calZ_a   \right \} .
\ee
  Clearly, this problem is equivalent to $(\calP_{sl})$ under the regularity condition, i.e.   $\bu = \bK(\brho)^{-1} \bff $ is well-defined on $\calZ_a$.
Instead, the given force is replaced by $\bff = \bK(\brho)\bu$ such that  $(\calP_c)$  is  commonly written  in the minimum strain energy form
\eb
 (\calP_s): \;\;\; \min
\left\{ \half \bu^T \bK( \brho) \bu \;\;  |  \;\;   \bK(\brho) \bu =  \bff ,   \;\;  \brho \in \calZ_a  \right\} .  \label{eq-mse}
\ee
Clearly, this problem  contradicts   $(\calP_{sl})$ in the sense that the alternative iteration for solving $(\calP_s)$
 leads  to an anti-Knapsack problem
\eb
 (\calP_{ak}): \;\; \min_{\brho}
\left\{   \bw^T  \brho|   \;\; \bv^T \brho \le V_c ,   \;\; \brho \in  \{0, 1\}^n \right\}. \label{eq-anti}
\ee
 By the fact that $\bw  = \bc(\bu)
=\{ \half   \bu_e ^T \bK_e \bu_e\}  \ge  {\bf 0} \in \real^n  \;\; \forall \bu\in \calU_a$, this anti-knapsack  problem has only trivial solution.

The   compliance $\bC$ is a well-defined  concept in engineering mechanics,  which is  complementary to the stiffness  $\bK$
 in the sense of $\bC= \bK^{-1}$.
In continuum physics, the linear scalar-valued function  $\bu^T \bff \in \real $ is called
the external (or input) energy, which is not  an objective function (see Appendix).
Since $\bff$ is a given force, it can't be replaced by $\bK(\brho) \bu$.
Although 
  the cost  function   $P_c(\brho) =  \half \bff^T \bC(\brho) \bff$ can be called  as the mean compliance, it is   not an objective function either. Thus,  the problem  $(\calP_c)$ works only for linear elastic structures. 
Its complementary form
\eb
 (\calP^c ): \;\;\; \max
\left\{ \half \bu^T \bK(\brho) \bu   \;\;  |  \;\;   \bK(\brho) \bu =  \bff ,     \;\;  \brho \in \calZ_a  \right\} \label{eq-msp}
\ee
 can be called a maximum stiffness problem, which  is equivalent to    $(\calP_{le})$ in the sense that both problems   produce the same results by  the alternative iteration method.
Therefore, it is a conceptual mistake to call  the strain energy $\half \bu^T \bK(\brho) \bu $ as the mean compliance and
  $(\calP_s)$  as the compliance minimization.\footnote{Due to this conceptual mistake,  the general problem for topology optimization was originally formulated as a double-min optimization
in \cite{g-to}. Although this model is equivalent to a knapsack problem for  linear elastic structures under the condition $\bff = \bK(\brho) \bu$, it contradicts the popular theory in topology optimization.}
Also,  the   compliance  can be defined only for linear elasticity.
 For nonlinear elasticity or plasticity, even if the stiffness can be defined as the Hessian matrix  $\bK(\bF) = \nabla^2 \WW(\bF)$, the associated compliance $\bC$ can't be well-defined since
$\bK(\bF)$ is usually not invertible due to the nonlinearity/nonconvexity  of the strain energy $\WW(\bF)$.

The   problem $(\calP_s)$ has been used extensively
by many  well-known  methods in topology optimization, including the most  popular SIMP
  \cite{Andreassen2011,sig-mau,bb_Zhou1991}:
 \eb
(\calP_{simp}): \;\;\;  \min_{\brho  \in (0, 1]^n  }
\left\{  \half \bu^T \bK( \brho^p) \bu  |  \;\;   \bK(\brho^p) \bu =  \bff ,  \;\; \bu \in \calU_a  ,  \;\;   \brho^T \bv \le V_c\right\}  \label{eq-simp}
\ee
Clearly,  the  integer constraint $\brho  \in \{ 0, 1\}^n$  in  $(\calP_s)$  is artificially replaced
 by a box constrain  $\brho  \in (0, 1]^n$ in  $(\calP_{simp})$  via  the so-called  power-law $\brho^p  = \{ \rho_i^p \} $.
 Although   $p >1$ is called the penalization parameter, the SIMP  is not a mathematically  correct penalty method for solving either  $(\calP_c)$ or $(\calP_s)$.
 In order to avoid the embarrassed anti-knapsack problem, the alternative iteration is not allowed by the SIMP method.
Therefore,  using  $\bu  =  \bC(\brho^p)   \bff$
the  strain energy  in $(\calP_{simp})$ is  written  as $P_s(\brho^p) = \half \bff^T  \bC(\brho^p)   \bff$.
Since $P_s(\brho^p)$ is not   coercive  on its domain,  unless some artificial techniques are adopted, the  global minimum solution of $(\calP_{simp})$
  can be achieved  only on its  boundary but can never  be   $0$   due to the restriction $\brho> {\bf 0} $.
The  so-called ``magic number'' $p=3$ works only for certain materials.
 This is the reason why the SIMP method suffers from  the  fatal drawbacks of
 gray scale  elements and checkerboard patterns.  \hfill $\clubsuit$

}
\end{remark}

\begin{figure} [h]
  \centering
  \includegraphics[width=0.9\textwidth]{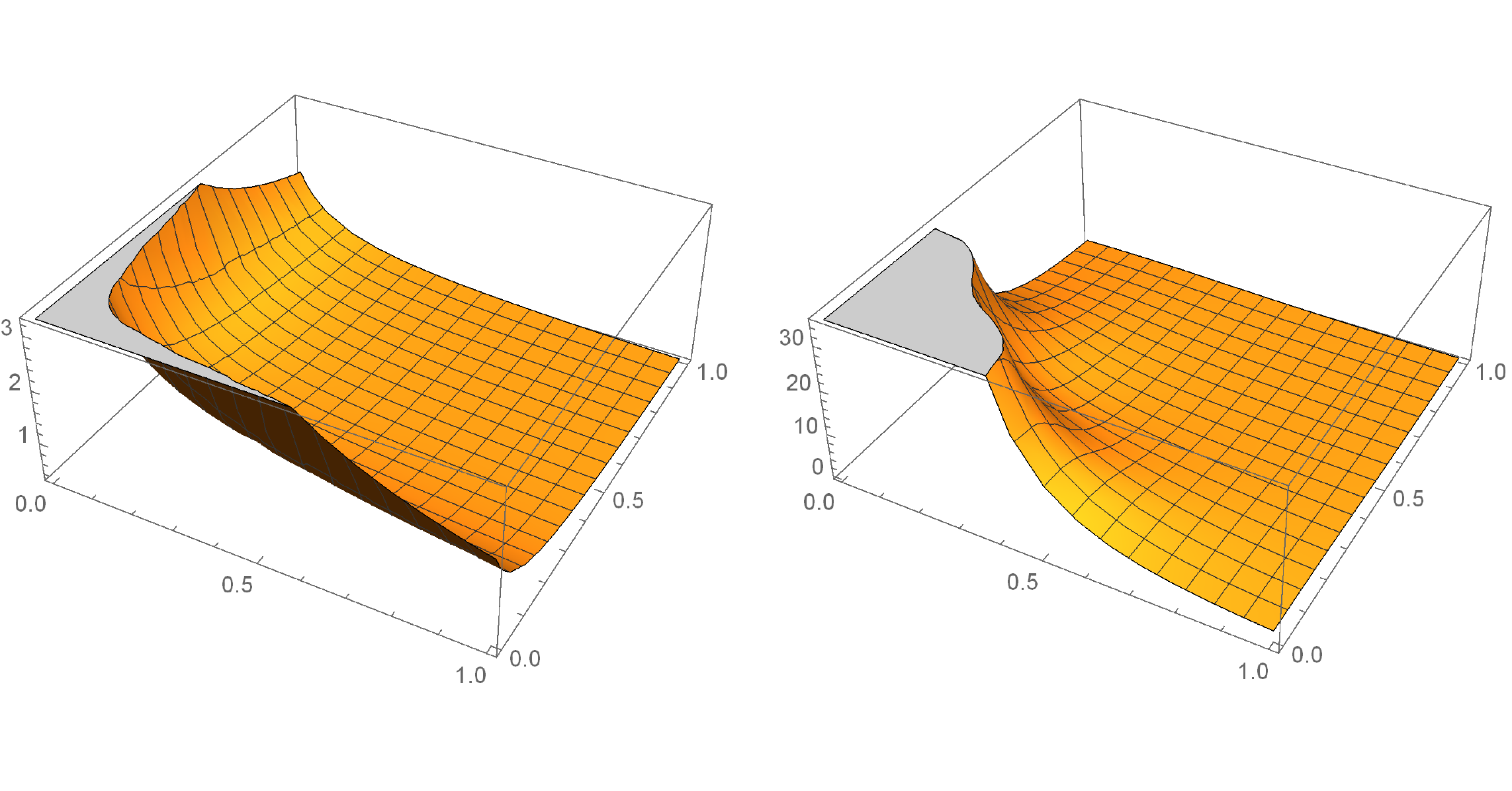}
  \caption{Graphs of $P_s(\brho^p)$ for $p=1$ (left) and $p=3$ (right)}
  \label{fig-ps}
\end{figure}

To understand this remark, let us consider a  2-D problem (see \cite{sto-sva} and Example  2 in Section 5) with
\eb
 \bK(\brho)   =  \left[ \begin{array}{cc}
a \rho_1 + b \rho_2 & 0 \\
0 & b \rho_1 + a \rho_2
\end{array} \right]
\ee
where $a, b \ge 0$ are material constants. For a given $\bff = \{f_i\}  \in \real^2$, the so-called penalized mean compliance is
\eb
P_s(\brho^p) =\half \bff^T \bK(\brho^p)^{-1} \bff =  \half [ {f_1^2}(a \rho_1^p + b \rho_2^p)^{-1} + {f_2^2 } (b \rho_1^p + a \rho_2^p)^{-1}  ]
\ee
Let $a =  (2 - \sqrt{2})/2, \; b = (4+ \sqrt{2})/2$, $f_1=f_2 = 1$, Fig. \ref{fig-ps} shows that $P_s(\brho^p) $ is not a coercive function  for $p=1,3$, and its global min can be obtained only at the boundary $\rho_1 + \rho_2 = 1$.
\begin{figure} [h]
 \centering
  \includegraphics[width=0.9\textwidth]{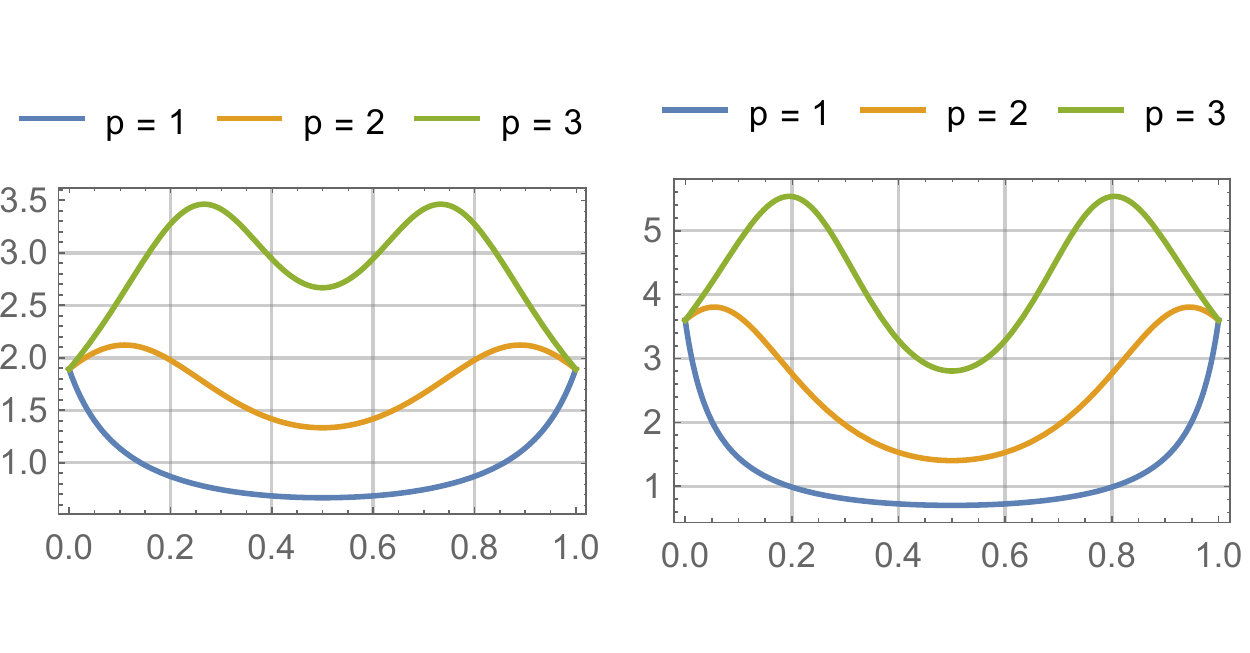}
\setlength{\unitlength}{.4cm}
\begin{picture}(-1,5)
\put(-30,-0){(a) $a =  (2 - \sqrt{2})/2$. \hspace{2cm} (b) $a =  (2 - \sqrt{2})/4$ }
 \end{picture}
  \caption{Graphs of $P_s(\brho^p)  $ on the boundary $\rho_1 + \rho_2 = 1$.}
  \label{fig-psb}
\end{figure}
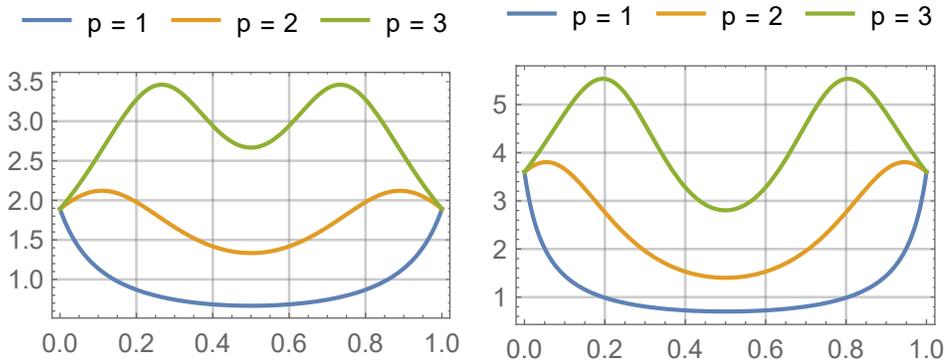
Fig. \ref{fig-psb} shows clearly that   $P_c(\brho ) $ is  strictly convex  on its boundary $\rho_1 + \rho_2 = 1$.
The feasible solutions  for   $(\calP_c)$ can be achieved  only on the
corners  $\brho = (1,0)$ and  $\brho=(0,1)$, both of them are global  maximizers. Therefore,  the problem
$(\calP_s)$   is considered to be NP-hard by  traditional   theories.
For the SIMP problem $(\calP_{simp})$ with $p=2$,
 its global minimizer  is  $ \brho=(0.5, 0.5)$, which is not a  solution to the topology optimization problem.
Although for  $p =3$ the global min of $P_s(\brho^p ) $ can be achieved at either
$\brho=(1,0)$ or  $\brho=(0,1)$ (see Fig. \ref{fig-psb} (a)), these two solutions can't be obtained by any Newton-type  method
since they are not a critical points.
Moreover,  if the  constant $a =  (2 - \sqrt{2})/4$, the global min of $P_s(\brho^3 ) $ is again $ \brho=(0.5, 0.5)$ (see Fig. \ref{fig-psb}(b)).
This shows that  $p=3$ in the SIMP method is a   ``magic number"   only for certain materials/structures.
By the fact that  $P_s(\brho^p ) $ is nonconvex for any given $p > 1$,
 the SIMP problem $(\calP_{simp})$
 is a typical ``NP-hard" box constrained  nonconvex minimization subjected to the knapsack condition,
which can't be solved deterministically by all  traditional methods in polynomial time.

\section{Canonical Dual  Solution to Knapsack Problem}\label{sec-dual}
 The canonical duality theory for solving constrained quadratic  0-1 integer programming problems was first proposed  by Gao in 2007~\cite{gao-jimo07}.
 Applications  have been given to general   problems in operations research \cite{gao-cace} and
 recently to topology optimization for general materials \cite{g-to}.
 In this paper, we   focus  mainly on linear elastic structures.

Following the standard procedure in the canonical duality theory (see~\cite{gao-jimo07,gao-ruan-jogo10}),
 the canonical measure for  0-1 integer programming problem $(\calP_{kp})$ can be given as
\eb
\bxi=  \Lam(\brho) = \{
\brho \ot  \brho - \brho, \;\; \brho^T \bv -V_c \} : \;\; \real^n \rightarrow
  \calE =  \real^{n+1}.
  \ee
  where $\brho \ot \brho= \{ \rho_e  \rho_e \} \in \real^n $ represents the Hadamard product.
Let
\eb
  \calE_a := \{ \bxi= \{ \beps, \nu\} \in \real^{n+1}|\;\; \beps \le 0, \;\; \nu \le 0 \}
\ee
be a convex cone in $\real^{n+1}$. Its indicator $\Psi(\bxi)$  is defined by
\begin{equation}\label{eq-ind}
\Psi(\bxi) = \left\{ \begin{array}{ll}
0 & \mbox{ if } \bxi \in \calE_a \\
+\infty & \mbox{ otherwise}
\end{array}
\right.
\end{equation}
which is a convex and lower semi-continuous (l.s.c) function in $\real^{n+1}$. By this function, the knapsack  problem
$(\calP_{kp})$ can be relaxed in the following unconstrained minimization form:
%
\eb\label{eq-problem-un}
\min  \left\{ \Pi_u (\brho) =   \Psi(\Lam(\brho)) -  \bw^T   \brho  \;\; | \;\;  \brho   \in \real^n  \right\}.
\ee

 By the convexity of   $\Psi(\bxi)$, its conjugate function can be defined uniquely  via the Fenchel transformation:
\eb\label{eq-psis}
\Psi^*(\bzeta) = \sup_{\bxi \in \real^{n+1}} \{ \bxi^T \bzeta - \Psi(\bxi) \}
=\left\{ \begin{array}{ll}
0 & \mbox{ if } \bzeta \in \calE_a^*\\
+\infty & \mbox{ otherwise},
\end{array} \right.
\ee
where
\eb\label{eq-as}
\calE_a^* = \{ \bzeta = \{ \bsig, \tau \} \in \real^{n+1} | \;\; \bsig \ge 0, \;\; \tau \ge 0 \}
\ee
is the  dual space of $\calE_a$.

Thus, by using the Fenchel-Young equality
\eb
\Psi(\bxi) + \Psi^*(\bzeta) = \bxi^T \bzeta,
\ee
the function $\Psi(\brho)$ in~\eqref{eq-problem-un} can be written as the Gao-Strang total complementary function~\cite{gs-89},
\eb\label{eq-gst}
\Xi(\brho, \bzeta) =     \Lam(\brho)^T \bzeta -  \bw^T   \brho - \Psi^*(\bzeta).
\ee

Based on this function, the \emph{canonical dual} of $\Pi_u(\brho)$ can be defined by
\eb\label{eq-cd}
 \Pi_u^d(\bzeta) = \mbox{sta } \{ \Xi(\brho, \bzeta) | \;\; \brho \in \real^n  \} = P^d_u (\bzeta ) - \Psi^*(\bzeta),
\ee
where
\eb
P^d_u (\bzeta ) = \mbox{sta } \{ \Lam(\brho)^T \bzeta  -  \bw^T   \brho | \;\; \brho \in \real^n \} =
 - \frac{1}{4} \btau^T_u(\bzeta ) \bG^{-1}(\bzeta)  \btau_u (\bzeta) - \tau V_c ,
\ee
  in which,
\[
\bG (\bzeta) = \ \Diag\{ \bsig\}  , \;\;\;  \btau_u ( \bzeta ) =    \bsig - \tau \bv  +  \bw .
\]

Clearly, $P_u^d(\bzeta)$ is well-defined if   $\bsig \neq {\bf 0 } \in \real^n $.
By the  standard  complementary-dual principle in the canonical duality theory, we have the  main result:
\begin{thm} [Complementary-Dual Principle] \label{thm-cdp}
For any given  $\bu \in \calU_a$, if $(\barbrho, \barbzeta) $ is a KKT point of $\Xi$, then $\barbrho$ is a KKT point of $ \Pi_u$, $\barbzeta$ is a KKT point of $ \Pi_u^d$, and
\eb
 \Pi_u (\barbrho) = \Xi(\barbrho, \barbzeta)= \Pi_u^d(\barbzeta).\label{eq-cdp}
\ee
\end{thm}
\emph{Proof:}
By the convexity of  $\Psi(\bxi)$, we have  the following canonical duality relations \cite{gao-jogo00,gao-dual00}:
 \eb \label{eq-cdr}
 \bzeta \in \partial \Psi(\bxi) \;\; \Leftrightarrow \;\; \bxi \in \partial \Psi^*(\bzeta)
 \;\; \Leftrightarrow \;\; \Psi(\bxi) + \Psi^*(\bzeta) = \bxi^T \bzeta,
 \ee
 where
 \[
 \partial \Psi(\bxi) = \left\{ \begin{array}{ll}
 \bzeta & \mbox{ if } \bzeta \in\calE_a^* \\
 \emptyset & \mbox{ otherwise}
 \end{array} \right.
 \]
 is the sub-differential of $\Psi$. Thus, in terms of $\bxi = \Lam(\brho)$ and $\bzeta = \{ \bsig, \tau\}$, the canonical duality relations (\ref{eq-cdr}) can be  equivalently written as
 \eb
 \brho \ot \brho - \brho \le 0   \;\; \Leftrightarrow \;\; \bsig \ge 0 \;\;
 \Leftrightarrow \;\; \bsig^T ( \brho \ot \brho - \brho )  = 0\;\; \label{eq-kkts}
 \ee
 \eb
 \brho^T \bv - V_c \le 0 \;\; \Leftrightarrow \;\; \tau \ge 0 \;\;\Leftrightarrow \;\;
 \tau(\brho^T \bv - V_c) = 0. \label{eq-kktv}
 \ee
These are exactly the KKT conditions\footnote{A   critical point  is a special  KKT point for equality constraint since the complementarity condition is automatically satisfied for all non zero Lagrange multipliers \cite{l-g-opl}.}
 for the inequality constraints $\brho \ot \brho - \brho \le 0  $ and $ \brho^T \bv - V_c \le 0$. Thus,   $(\barbrho, \barbzeta) $ is a KKT point of $\Xi$ if and only if
$\barbrho$ is a KKT point of $ \Pi_u$, $\barbzeta$ is a KKT point of $ \Pi_u^d$. The equality (\ref{eq-cdp}) holds due to the canonical duality relations in (\ref{eq-cdr}). \hfill  \qed


By the complementarity condition $\bsig^T ( \brho \ot \brho - \brho ) = 0$ in \eqref{eq-kkts}, we know that $\brho \ot \brho = \brho  $  if  $\bsig > 0$.
Let
\eb
\calS_a^+ = \{ \bzeta = \{ \bsig, \tau\}\in  \calE_a^* | \;\;    \bsig  >  0 \}.
\ee
Then for any given  $ \bzeta= \{ \bsig, \tau\}  \in \calS^+_a$,  the function $\Xi(\cdot, \bzeta): \real^n \rightarrow \real$ is  convex and
the canonical dual function of $\PP_u$ in $(\calP_{kp})$ can be well-defined by
\eb\label{eq-Pd}
\PP^d_u(\bzeta) = \min_{\brho \in\real^n} \Xi(\brho, \bzeta) =
- \sum_{e=1}^n  \frac{1}{4} (\sig_e  + \ww_e
- \tau \vv_e)^2 \sig_e^{-1}  - \tau V_c.
\ee
Thus, the canonical dual problem of $(\calP_{kp})$ can be proposed as the following:
\eb
(\calP^d_u): \;\;\; \max \{ \PP^d_u(\bsig,\tau) | \;\; (\bsig, \tau) \in \calS^+_a \}.
\ee

This is a concave maximization problem over a convex subset in continuous space, which can be solved via well-developed convex optimization  methods  to obtain global optimal solution.
 Thus,  whence a canonical dual solution $  \barbzeta $ is obtained,
 the  solution  $\bar{\brho}$ to the primal problem  can be defined  in  an  analytical form of $  \barbzeta $.
\begin{thm}[Analytical Solution] \label{thm-rho}
For any given $\bu \in \calU_a$ such that $\bw = \bc(\bu) \in \real^n_+$, if $  \barbzeta = (\barbsig, \bartau) \in \calS^+_a  $ is a  solution to $(\calP^d_u)$,
 then
\eb\label{eq-solu}
\barbrho =
\half [\Diag(\barbsig)]^{-1} (\barbsig -\bartau \bv +  \bw  )
\ee
is a unique global optimal solution to $(\calP_{kp})$ and
\eb
\PP_u(\barbrho) = \min_{\brho \in \calZ_a} \PP_u(\brho) = \max_{\bzeta \in \calS^+_a } \PP^d_u(\bzeta)
=  \PP^d_u(\barbzeta).
\ee
\end{thm}
\emph{Proof}:
 It is easy to prove that   for any given  $\bu \in \calU_a$, the canonical dual function $\PP^d_u(\bzeta)$ is concave on the open convex set $\calS^+_a$. If $\barbzeta$ is a KKT point of $\PP^d_u(\bzeta)$, then it must be a unique global maximizer of $\PP^d_u(\bzeta)$ on $\calS^+_a$.
  By Theorem 1 we know that if $\barbzeta = \{ \barbsig , \barvsig \} \in \calS^+_a$ is a KKT point of $\Pi_u^d(\bzeta)$, then $\barbrho = \brho(\barbzeta)$ defined by (\ref{eq-solu}) must be a KKT point of $\Pi_u(\brho)$ (see \cite{gao-jimo07,g-to}). Since $\Xi(\brho, \bzeta)$ is a saddle function on $\real^n \times \calS^+_a$, we have
\begin{eqnarray*}
 \min_{\brho \in \calZ_a } \PP_u(\brho)  &=& \min_{\brho \in \real^n}  \Pi_u (\brho)=\min_{\brho \in \real^n} \max_{\bzeta \in \calS^+_a}
\Xi(\brho, \bzeta) =   \max_{\bzeta \in \calS^+_a}\min_{\brho \in \real^n}
\Xi(\brho, \bzeta)\\
& =& \max_{\bzeta \in \calS^+_a}  \Pi_u^d(\bzeta) =  \max_{\bzeta \in \calS^+_a} \PP^d_u(\bzeta).
\end{eqnarray*}
Since $\barbsig > 0$,  the complementarity condition in (\ref{eq-kkts}) leads to
 \[
 \barbrho \ot \barbrho - \barbrho = 0 \;\;  \mbox{ i.e. } \barbrho \in \{ 0, 1 \}^n.
 \]
 Thus, we have
 \[
 \PP_u(\barbrho) = \min_{\brho \in \calZ_a } \PP_u(\brho) = \max_{\bzeta \in \calS^+_a } \PP^d_u(\bzeta)
 = \PP^d_u(\barbzeta)
 \]
 as required. \qed

Theorem \ref{thm-rho} shows that although the canonical dual problem is a concave maximization in continuous space,  it produces the analytical solution  (\ref{eq-solu})
to  the well-known integer Knapsack problem $(\calP_{kp})$! This analytical solution was first proposed by Gao in 2007
for  general quadratic integer programming problems (see Theorem 3, \cite{gao-jimo07}).
The indicator function of a convex set in~\eqref{eq-ind} and its sub-differential were first introduced by J.J. Moreau in 1963 in his study on unilateral constrained problems in contact mechanics \cite{moreau}. His pioneering work laid a foundation for modern analysis and the canonical duality theory.  In solid mechanics, the indicator of a plastic yield condition is also called a {\em super-potential}. Its sub-differential leads to a general constitutive law and a unified pan-penalty finite element method in plastic limit analysis \cite{gao-cs88}. In mathematical programming, the canonical duality leads to a unified framework for nonlinear constrained optimization problems in multi-scale systems \cite{gao-ruan-sherali-jogo,l-g-opl}.

\section{Canonical Penalty-Duality Method and  Algorithm}\label{sec-num}
According to Theorem~\ref{thm-rho}, the global optimal solution  to $(\calP_{kp})$ can be obtained by solving
 its canonical  dual problem $(\calP_{kp}^d)$. However, the rate of convergence could be
  very slow since $\PP^d_u(\bsig,\tau)$ is nearly a linear function of $\bsig \in \real^n_+$ when $\bsig $ is far from its origin. In order to overcome this problem,
 a so-called $\beta$-perturbed canonical dual method was  proposed by Gao and Ruan in integer programming \cite{gao-ruan-jogo10}.
This  $\beta$-perturbation is actually based on a so-called canonical penalty-duality method, i.e.  the integer constraint in $(\calP_{kp})$
 is first written in the canonical form $\beps = \brho \ot\brho - \brho = {\bf 0} $, then  is relaxed by the external penalty method, while the volume constraint is simply relaxed by the Lagrange multiplier $\tau \ge 0$, thus, the knapsack problem can be reformulated as the following canonical penalty-duality form
 \eb
\min_{\brho  \in \real^n} \max_{\tau \ge 0 } \left\{ P_\beta(\brho, \tau) = \beta \| \brho \ot \brho - \brho \|^2 - \bw^T \brho + \tau (\bv^T \brho - V_c) \right\},
\ee
where $\beta > 0$ is a penalty  parameter.
As we can see clearly that due to the nonlinearity of the canonical  constraint $\beps ( \brho)  = 0$, the canonical penalty function  is nonconvex in $\brho$. This is  the reason
why  the  standard penalty method can't be used for solving general nonlinearly constrained problems.
However, by using the canonical duality theory, this nonconvex min-max  problem can be equivalently reformulated as a canonical dual problem:
\eb
 (\calP^d_\beta): \;\; \max \left\{ \PP^d_\beta(\bsig, \tau)  = \PP^d_u (\bsig, \tau) -
\frac{1}{4} \beta^{-1} \bsig^T  \bsig | \;\; \{ \bsig, \tau\} \in \calS^+_a \; \right\} ,
\ee
which is exactly the same canonical penalty-duality problem  proposed by Gao for solving the knapsack problem in topology optimization  \cite{g-to}.
\begin{thm}[Perturbed  Solution to Knapsack Problem] \label{thm-solution} For any  given $\bu \in \calU_a$, $V_c > 0$, and $\beta  > 0$ such that  $\bw=\bc(\bu) \in \real^n_+$,
the problem $ (\calP^d_\beta))$ has at most one solution
  $\bzeta_\beta = (\bsig_\beta, \tau_\beta)  \in \calS^+_a$. Moreover, there exists a $\beta_c  \gg 0 $ such that for any given $\beta \in [\beta_c, +\infty)$,  the vector
\eb\label{eq-rhob}
\brho_\beta = \half [\Diag(\bsig_\beta) ]^{-1} (\bsig_\beta  - \tau_\beta \bv + \bw)  \in\{ 0, 1\}^n
\ee
is   a global optimal solution to $(\calP_{kp})$.
\end{thm}

\emph{Proof}. It is easy to show that for any given $\beta > 0$,
\eb
\PP^d_\beta(\bsig, \tau) = - \frac{1}{4} \left [ \sum_{e=1}^n  (\sig_e  + \ww_e
- \tau \vv_e)^2 \sig_e^{-1}  +  \beta^{-1} \sig^2_e \right] - \tau V_c
\ee
is strictly concave on  the open convex set $\calS^+_a$. Thus,  $(\calP^d_\beta)$ has a unique solution
 if $\ww_e
- \tau \vv_e \neq 0 \;\; \forall e=1, \dots, n$.
Indeed, the criticality condition $\nabla \PP^d_\beta(\bzeta) = 0 $ leads to the following canonical dual algebraic equations:
\eb
2 \beta^{-1} \sig_e^3 + \sig_e^2 = (\tau \aa_e - \ww_e )^2, \;\; e = 1, \dots, n, \label{eq-cdas}
\ee
\eb
\sum_{e=1}^n \half \frac{\aa_e}{\sig_e} ( \sig_e - \aa_e \tau + \ww_e ) - V_c = 0 .\label{eq-cdv}
\ee
It was  proved by the author (see Section 3.4.3, \cite{gao-dual00} and the Appendix of this paper)  that for any given $\beta > 0$
and $\theta_e (\tau) = \tau \aa_e - \ww_e   \neq 0, \ e = 1, \dots, n$, the canonical dual algebraic equation (\ref{eq-cdas}) has a unique
positive real solution
\eb
\sigma_e  =  \frac{1}{6} \beta   [- 1 +  \phi_e(\tau  ) + \phi_e^c(\tau  )] > 0 , \;\; e = 1, \dots, n
\label{eq-solus}
\ee
where
\[
\phi_e(\tau )  = \eta^{-1/3} \left [2 \theta_e(\tau) ^2 - \eta + 2 i     \sqrt{\theta_e (\tau)^2 (\eta -  \theta_e(\tau) ^2 )
 } \right]^{1/3} ,
 \;\; \eta = \frac{\beta^2}{27 },
\]
and $\phi_e^c $ is the complex conjugate of $\phi_e $, i.e. $\phi_e  \phi_e^c  = 1$.
Also,  the canonical dual algebraic equation (\ref{eq-cdv}) has a unique solution
\eb\label{eq-soluvs}
 \tau = \frac{\sum_{e= 1}^n \aa_e ( 1 +   \ww_e /\sig_e) - 2 V_c}{\sum_{e=1}^n \aa_e^2/\sig_e}.
\ee
This shows that the perturbed canonical dual problem $(\calP^d_\beta)$ has a unique solution $\bvsig_{\beta}$ in $\calS^+_a$ if
 $\theta_e (\tau) = \tau \aa_e - \ww_e  \neq 0, \ e = 1, \dots, n$.
  Thus, the density distribution $\brho_{\beta}$
can be analytically obtained by substituting~\eqref{eq-solus} and~\eqref{eq-soluvs} into~\eqref{eq-rhob}.
By the fact that $\lim_{\beta \rightarrow \infty } \Pi^d_\beta(\bzeta) = \Pi^d_u(\bzeta)$, there must exists a $\beta_c > 0$ such that
  \eb
  \brho_\beta  \in\{ 0, 1\}^n  \;\; \forall \beta \ge \beta_c .
  \ee
  By Theorem \ref{thm-rho} we know that this perturbed solution must be a global minimum solution to the knapsack problem.
A similar proof of this theorem was  given in \cite{gao-ruan-jogo10}.
\hfill $\Box$


By  Theorems  \ref{thm-rho} and \ref{thm-solution}  we know that for a  given desired volume $V_c> 0$, the optimal density distribution  can be analytically obtained
 in terms of its canonical dual solution in continuous space. By the fact that   $(\calP_{bk})$ is a
 bilevel mixed integer nonlinear programming, numerical optimization  depends   sensitively on
 the  initial volume $V_0$. If $ \mu_c = V_c/V_0 \ll 1, $ any given iteration  method  could lead to unreasonable numerical solutions. In order to resolve this problem,
a volume reduction  control parameter
$\mu \in (\mu_c,1)$ was introduced in \cite{g-to} to produce a volume sequence  $V_{\gamma } = \mu V_{\gamma-1}$ ($\gamma =  1, \dots, \gamma_c$)
such that $V_{\gamma_c} = V_c$ and for any given $V_\gamma \in [V_c, V_0]$, the problem $(\calP_{bl})$ is replaced by
\begin{eqnarray}
 (\calP_{bk})^\gamma:\;\; &    &  \min
\{    \bff^T \barbu  - \bc(\barbu )^T \brho^p |\;\; \bv^T \brho \le V_\gamma, \;\;  \brho \in \{ 0, 1\}^n  \}   \\
 & & \mbox{s.t.}  \;\; \barbu(\brho)  =\arg \min \{  \Pi_h(\bu, \brho) | \;  \;{\u \in \calU_a }\}.
\end{eqnarray}
 The initial values for solving this $\gamma$-th problem are $V_{\gamma-1}, \bu_{\gamma-1}, $ $ \brho_{\gamma-1}$.
Based on the above strategies, the canonical duality algorithm (i.e. {CDT} \cite{g-to}) for solving the general topology optimization problem $(\calP_{bk})$ can be proposed in Algorithm 1.

\begin{algorithm}
\caption{Canonical Dual Algorithm for Topology Optimization (CDT)}
\begin{algorithmic} [1]
\State  \textbf{Input parameters:}  $\mu$,   $\beta$ and error allowances: $\omega_1,  \omega_2 > 0$.
\State  \textbf{Initiate} $\brho^0 = \{1\} \in \real^n$,   $V_0 , \tau^0 > 0$, and   let $\gamma = 1 $.
\State  \textbf{Solve}  the lower-level optimization problem  \eqref{eq-lowk}:
    \[
      \u^{\gamma }   = \arg \min\{ \Pi_h(\bu, \brho^{\gamma-1})  | \;\; \bu \in \calU_a \}.
    \]
\State \textbf{Compute} $\bw^{\gamma } = \{ \ww_e ^\gamma \} = \bc(\bu^{\gamma })$  via \eqref{eq-cu}, $V_{\gamma } = \max\{ V_c,  \mu V_{\gamma-1}\} $.   Let $ k=1$.
  \State \textbf{Compute the canonical dual solutions}   $\bsig^k = \{\sigma_e^{k } \}$ and $    \tau^{k } $ by
        \[
        \sigma_e^{k } =  \frac{1}{12} \beta   [- 1 +  \phi(\tau^{k-1}  ) + \phi^c(\tau^{k-1 } )] , \;\; e = 1, \dots, n.
        \]
        \[
        \tau^{k } = \frac{ \sum_{e=1}^n  \aa_e (1 + \cc^\gamma_e /\sigma_e^{k} ) - 2 V_{\gamma } }{\sum_{e=1}^n \aa_e^2/  \sigma_e^{k} }.
        \]
   \State {\bf  If }  $|P^d_u(\bsig^k, \tau^k)  - P^d_u( \bsig^{k-1},\tau^{k-1}) |  >  \omega_1$ , then   let  $k = k + 1$, go to Step 5;
    Otherwise,   continue.
     \State   \textbf{Compute the upper-level  solution} $\brho^{\gamma } $ by 
        \[
        \rho^{\gamma }_e   = \frac{1}{2} [ 1 - ( \tau^{k } \aa_e - \cc^\gamma_e)/\sigma_e^k],
        \;\; e= 1, \dots, n.
        \]
      \State    {\bf  If }   $| P_u( \brho^{\gamma })  -  P_u(\brho^{\gamma-1})   | \le \omega_2$ and $V_\gamma  \le V_c$ ,
      then stop; Otherwise, continue.
                \State \textbf{Let}    $\tau^0 = \tau^{k }$,   and  $\gamma = \gamma+1 $. Go to Step 3.
\end{algorithmic} \label{alg-CDT}
\end{algorithm}

\begin{remark} [Volume Reduction and Computational Complexity] $\;$\newline
 {\em Theoretically speaking,  for any given sequence $\{V_\gamma \}$ we should have
 \eb
 (\calP_{bk}) = \lim_{\gamma\rightarrow \gamma_c} (\calP_{bk})^\gamma.
 \ee
In reality, different sequence $\{V_\gamma\}$   may produce totally different structural topology.
This is an intrinsic difficulty for all   bi-level  optimal design problems.
  By the facts  that   there are only two loops in the CDT  algorithm, i.e. the $\gamma$-loop and the $k$-loop,  and   the canonical dual solution is analytically given in the $k$-loop,
the main computing is the $m\times m$  matrix inversion in the $\gamma$-loop. The complexity for the Gauss-Jordan elimination is $O(m^3)$. Therefore, the CDT
  is a  polynomial-time algorithm.

The   optimization problem of BESO as formulated in
\cite{bb_Huang2010} is   posed in the form of   minimization of  mean
compliance, i.e. the problem $(\calP_c)$.
Since   the alternative iteration is adopted by BESO, which  leads to    an anti-Knapsack problem by \eqref{eq-anti} , therefore,
the BESO  should theoretically produce only trivial solution at each volume evolution.
However, instead of solving the anti-Knapsack problem (\ref{eq-anti}),
  a   comparison method is used   to  keep  those elements which store more strain energy. So, the BESO
 is  actually a direct method for solving
    the knapsack problem $(\calP_{kp})$.  This is the reason why the numerical results obtained by BESO are similar to that by CDT (see Section 6).
But,  the  direct  method  is not a  polynomial-time algorithm.  Due to the combinatorial complexity,  this popular method is computationally expensive and
can  be used only for small sized problems.
This is the very  reason why the knapsack problem has been  considered as  NP-complete for  all existing direct approaches. \hfill$\clubsuit$
}
 \end{remark}
\vspace{-.5cm}
\section{Symmetry and NP-Hardness for Knapsack Problem} \label{sec4}
There is a hidden  condition in   the proof of Theorem  \ref{thm-solution},  which is actually the limitation of the
canonical duality theory. 

\begin{thm}[Existence of Solution to Knapsack Problem] \label{thm4}
For any given   $\bv, \bw \in \real^n_+$, if there exists a constant $\tau_c >  0$ such that
\eb
\theta_e (\tau_c) = \tau_c  \vv_e  - \ww_e \neq { 0} \;\; \forall e=1, \dots, n \label{eq-exis}
\ee
the canonical dual  feasible set  $\calS^+_a \neq \emptyset $ and     the   knapsack problem $(\calP_{kp})$  has a  unique solution.
   Otherwise,  if $\theta_e(\tau_c) = 0 $ for at least one $e \in \{1, \dots, n\}$, then $\calS^+_a = \emptyset$ and $(\calP_{kp})$  has at least two  solutions.
   \end{thm}
   \emph{Proof}. From the proof of Theorem  \ref{thm-solution} we know that for any given $\beta> 0$, if
   thre exists a $\tau_c > 0$ such that the conditions in  \eqref{eq-exis}  hold,
    the canonical dual algebraic equation  \eqref{eq-cdas} has  a unique $\sig_e > 0$ for every $e = 1, \dots, n$,
    which can be obtained analytically by \eqref{eq-solus} see Theorem 3.4.4 in \cite{gao-dual00}).
    Therefore, the canonical dual problem $(\calP^d_{kp})$ has a unique solution in $\calS^+_a$.
    Correspondingly, the primal problem  $(\calP_{kp})$ has a unique solution defined by  either  \eqref{eq-solu} or \eqref{eq-rhob}.
   If  $\theta_e(\tau_c) = 0 $  for at least one $e = 1, \dots, n$,  the  equation 
 $ 2 \beta^{-1} \sig_e^3 + \sig_e^2 = 0$ has two solutions $\sig_e = 0$. In this case, the KKT points of
     $(\calP^d_\beta)$  are located on the boundary of the open set $\calS^+_a$. Therefore, $\calS^+_a = \emptyset$ and $(\calP_{kp}) $   has at least two solutions.
   \hfill $\Box$

    \begin{remark}[NP-Hard Conjecture, Symmetry, and Linear Perturbation]
  {\em   Theorem \ref{thm4}  shows that  under the condition \eqref{eq-exis}, the well-known knapsack problem is not NP-hard and can be solved analytically by the canonical duality theory.
It is discovered recently \cite{g-18} that the  critical value $\tau_c$ in  \eqref{eq-exis} can be deterministically  given by
\eb
\tau_c = \arg \min_{\tau \ge 0 }  \left\{  \sum_{e=1}^n ( | \ww_e - \tau \vv_e | - \tau \vv_e)  + 2 \tau   V_c  \right\} . \label{eq-tauc}
\ee
   Otherwise, as long as  $\theta_e(\tau_c) = 0 $  for at least one $e = 1, \dots, n$,
     the solution to $(\calP_{kp})$ can't be written in the analytical form   \eqref{eq-rhob},  and
  the knapsack problem could be really NP-hard.
  Actually, it is a conjecture first  proposed by the author  in 2007 \cite{gao-jimo07}, i.e.
  \begin{verse}
{{\bf Conjecture of NP-Hardness}:   A  global optimization  problem is NP-hard if its canonical dual problem has no solution in $\calS^+_a$. }
\end{verse}
  It is also an open problem left in \cite{gao-cace,gao-ruan-jogo10}.
The reason for   NP-hard problems and possible solutions  were discussed recently in \cite{gao-aip}.

 Geometrically speaking, the reason for  multiple solutions of $(\calP_{kp})$  is due to certain symmetry on the modeling, boundary condition and the  external load.
 By the fact that nothing is perfect in this real world, a perfect symmetry is not allowed for any real-world problem.
Mathematically speaking,  if   a problem  has multiple solutions, this problem is not {\em well-posed} \cite{gao-opl16}.
 In order to solve such NP-hard problems, the key idea is to break the symmetry.
 A linear perturbation method has been  proposed by the author and his co-workers with successful applications
 in hard cases of trust region method \cite{chen-gao-amma}, nonconvex constrained optimization \cite{mor-gao-amma}, and  integer programming \cite{wangetal}. \hfill$\clubsuit$
}
\end{remark}

 The symmetry plays a fundamental rule in mathematical modeling. But, it is also a main reason that leads to chaos in nonlinear dynamics, post-buckling in large deformation mechanics,  NP-hard problems in complex systems \cite{g-l-r-17,l-g-chaos}, and the well-known  paradox of
  Buridan's ass\footnote{Jean Buridan,  (born 1300, probably at B\'{e}thune, France--died 1358),
Aristotelian philosopher, logician, and scientific theorist in optics and mechanics.}. 

\begin{exam}[Buridan's Ass]
{\em  A donkey facing 
  two identical hay piles starves to death because reason provides no grounds for choosing to eat one rather than the other.
  Mathematically, this is   a knapsack problem:
  \eb
  \max \{  \ww_1 \rho_1 + \ww_2 \rho_2  |\;\;\; \ww_1 = \ww_2 = \ww, \;\; \rho_1 + \rho_2 \le 1, \;\; (\rho_1, \rho_2)  \in \{ 0,1\}^2 \}
  \ee
Due to the symmetries: $\vv_1 =\vv_2 = 1, $ and $ \ww_1=\ww_2 = \ww$, the solution to \eqref{eq-tauc}
 is $\tau_c = \ww$. Therefore, $\theta_e(\tau_c) = 0 \;\; \forall e = 1,2$ and by Theorem \ref{thm4}
this problem has multiple (two) solutions, which is NP-Hard to this donkey.

To solve this problem,  a linear perturbation term $\eps \rho_1  $  can be added to the cost function to break the symmetry. For $\ww=2, \;\; \eps=0.05$, we have $\tau_c = 2.0184$.
So the condition \eqref{eq-exis}  holds for $e=1,2$ and by the canonical duality theory,  the perturbed Buridan's ass problem  has a unique solution $\brho = (1,0)$.
}
\end{exam}

\begin{exam}[Symmetrical Truss Topology Optimization]
{\em
Let us consider a structure proposed in \cite{sto-sva} (see Fig. \ref{sto-sva})  with six bars that
 are grouped into 2 groups. Each group has the same cross sectional area, which is the design variable. The group   stiffness matrices  and load of this structure are:
 \begin{figure} [h]
  \centering
  \includegraphics[width=0.6\textwidth]{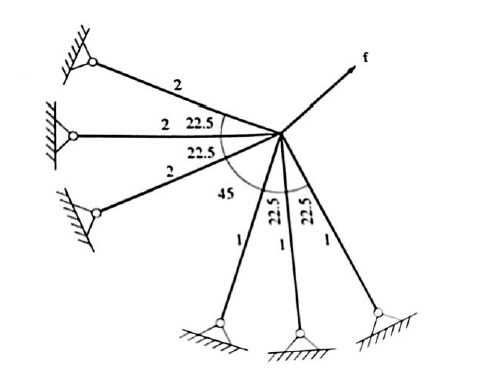}
  \caption{Symmetrical truss  with symmetrical load}
  \label{sto-sva}
\end{figure}\[
 \bK_1 = \left[ \begin{array}{cc} a & 0 \\ 0 & b \end{array} \right], \;\;
 \bK_2 =  \left[ \begin{array}{cc} b & 0 \\ 0  &  a\end{array} \right],
\;\;
\bff =  \left[ \begin{array}{c} 1   \\  1  \end{array} \right], \;\; a = \frac{2-\sqrt{2}}{2}, \;\;
b= \frac{4+ \sqrt{2}}{2}.
\]
For this linear elastic system, the associated  problem  $(\calP_{le})$ is
\begin{eqnarray}
 (\calP_{le}): \;\;  &  \min    &
\{  \bff^T \bu - \brho^T \bc( \bu)  | \;\;  \;\; \rho_1+ \rho_2 \le 1, \;\; \brho \in \{ 0, 1\}^2  \}   \label{eq-kna} \\
&  \mbox{s.t.} &  (  \rho_1   \bK_1 +   \rho_2   \bK_2 ) \bu = \bff  ,
\end{eqnarray}
where $ \bw = \bc(\bu) = \half ( a \uu_1^2  + b \uu_2^2 , \;\; b \uu_1^2  + a \uu_2^2 ) $ .
Due to the symmetry, for a given $V_0=2$ and  $\brho_0 = (1,1)$  we have $\bu_0 = ( 1/(a+b) , 1/(a+b) )$, $\bc(\bu_0) = ( \ww_1, \ww_2)$ with $\ww_1 =\ww_2 = 1/(a+b)$.
Thus,  the
upper-level optimization \eqref{eq-kna} for the first iteration  is exactly the Buridan  ass problem.
Clearly, this problem is artificial as it impossible to have a perfectly symmetrical truss with the perfect load $\bff = (1,1)$.
By using linear perturbation $\bff_\eps = (1+ \eps , 1)$, for any $\eps > 0$ we have $\uu_1 > \uu_2$ and $\ww_1 < \ww_2$. The CDT algorithm can produce a unique global optimal solution
$\brho^+ = (0,1)$. Dually, for $\eps < 0$, we have $\brho^- = (1,0)$, and we have
$\Pi_h(\bu(\brho^+) , \brho^+)= \Pi_h(\bu(\brho^-) ,\brho^-) = - \half (1/a + 1/b)$.
}\end{exam}


\section{Applications to Benchmark Problems and Novelty}\label{sec-example}
The proposed CDT algorithm  for topology optimization
 has been implemented in Matlab. The  CDT  code for  2-D topology optimization is based on the popular  88-line SIMP code  (TOP88) proposed by Andreassen et al~\cite{Andreassen2011}. By the facts  that the density distribution is solved analytically at each iteration and  no density filter is needed,   the CDT code has only 66-lines.
 The CDT code  for 3-D topology optimization  is based on the TOP3D  code  proposed  by Liu and Tovar \cite{liu-tovar}.
 The CDT code has been performed for  various numerical examples to test its performance.
 For the purpose of  illustration, the  applied  load and geometry data are chosen as dimensionless.
Young's modulus and Poisson's ratio of the  material are taken as $E = 1$ and $\nu = 0.3$, respectively.
The stiffness matrix of the structure in CDT algorithm  is given by
\[
\bK (\brho) = \sum_{e=1}^n [ E_{\min}  + (E-E_{\min} ) \rho_e ]\bK_e ,  \;\; E_{\min} = 10^{-9} .
\]
Clearly, we have
$\bK ({\bf 1}) = \sum_{e=1}^n   E \bK_e$  and $\bK ({\bf 0}) = \sum_{e=1}^n   E_{\min}  \bK_e$.
The reason for  choosing     $E_{\min} \neq 0 $ is  to avoid singularity in computation.
To compare with other approaches,  the
parameters  penal $=3$,
rmin = 1.5, and  ft=1.0  are used in
  the SIMP   88-line code, BESO code,  and  the  TOP3D  code.
 The error allowances are  set to be  $\omega_1 = 2e-16$ for CDT algorithm and
 $\omega_2 = 10^{-2}$ for all methods (SIMP is usally failed to converge if $\omega$ is too small). 
The initial value  for   $\tau$ used in CDT   is $\tau^0=1$.
We take the design domain $V_0 = 1$, the initial design variable  $\brho^0=\{1\}^n$ for both  CDT and BESO algorithms.
All computations are performed by a  HP laptop computer with Processor Intel Core I7-4810,  CPU @ 2.80GHz and memory 2.80 GB.

\subsection{MBB  Beam Problem}

\begin{figure} [h]
  \centering
  \includegraphics[width=0.6\textwidth]{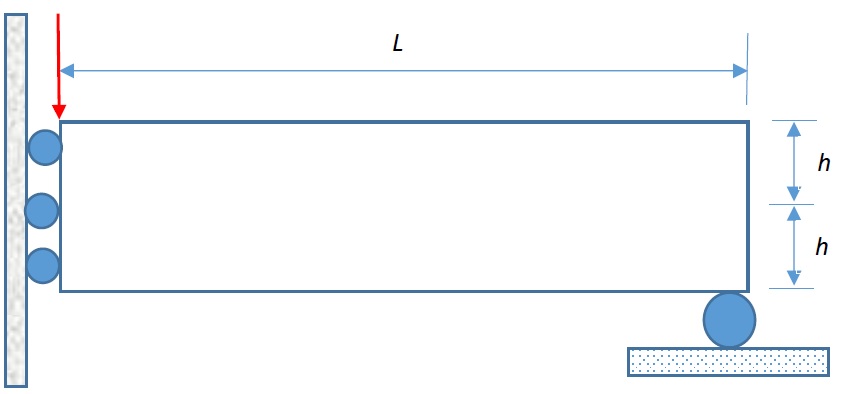}
  \caption{The design domain   for  a half MBB  beam with external load}
  \label{fig-mbb}
\end{figure}

\begin{figure} [h]
  \centering
  \includegraphics[width=0.6\textwidth]{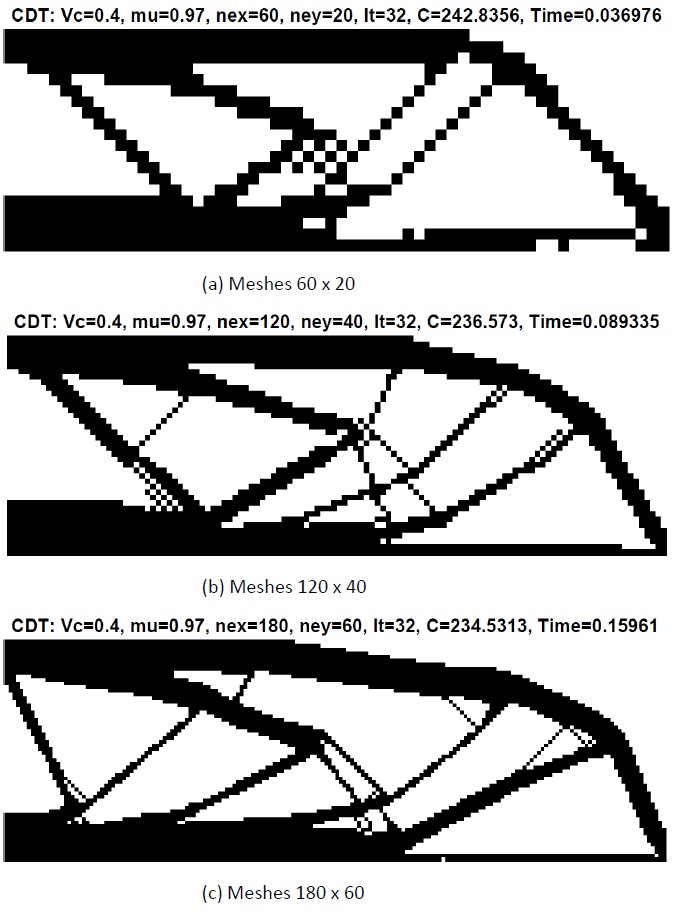}
  \caption{Optimal structures by different mesh resolutions with  $\mu_c = 0.4$ at  $\mu=0.97$.}
   \label{fig2}
\end{figure}

The first example is  the   well-known benchmark Messerschmitt-B\"{o}lkow-Blohm (MBB) beam problem in topology optimization
    (see  Fig.~\ref{fig-mbb}).
The design domain is  $L\times 2 h = 180 \times 60$. 
 Performance of the  CDT method  is first  tasted   for different mesh resolutions.
Results in Fig. \ref{fig2} show    that for any given mesh resolutions,
 the CDT method produces precise integer solutions without using  filter.
Clearly,   the finer the resolution, the smaller  the compliance with better result (almost no checkerboard).

\begin{figure} [h]
  \centering
  \includegraphics[width=1\textwidth]{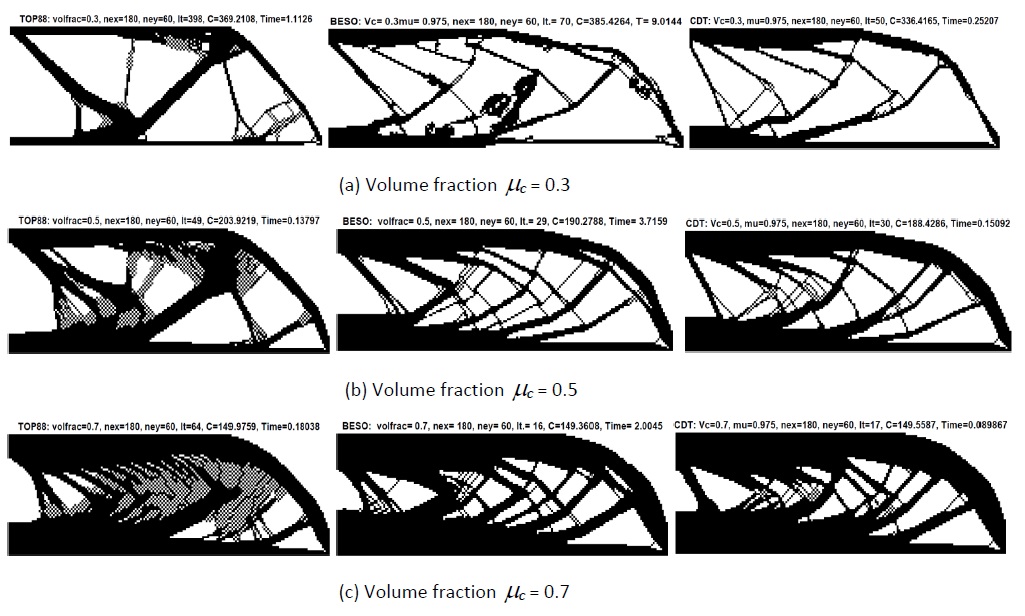}
  \caption{Computational results by SIMP (left), BESO (middle), and CDT (right) with different volume fractions $\mu_c = V_c/V_0$. For $\mu_c = 0.3$, the TOP88 code  failed to  converge. The result reported in (a) is the out put at the iteration It = 398.   }
  \label{fig3}
\end{figure}

To compare with the SIMP and BESO methods, we use the same mesh resolution of $180\times 60$ but with different volume fractions.
The volume reduction  rate is fixed to be $\mu=0.975$ for both CDT and BESO.
Computational results are reported in
Figure \ref{fig3}, which show clearly that for any given volume fraction $\mu_c$, the CDT method produces better results within the significantly short times.

 \subsection{2-D Cantilever Beam Problem}
\begin{figure} [h]
  \centering
  \includegraphics[width=0.6\textwidth]{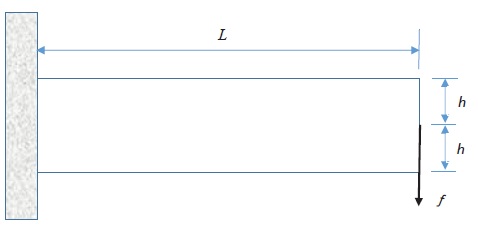}
  \caption{The design domain   for  a long cantilever   beam with external load}
  \label{fig-cant}
\end{figure}

\begin{figure} [h]
  \centering
  \includegraphics[width=0.6\textwidth]{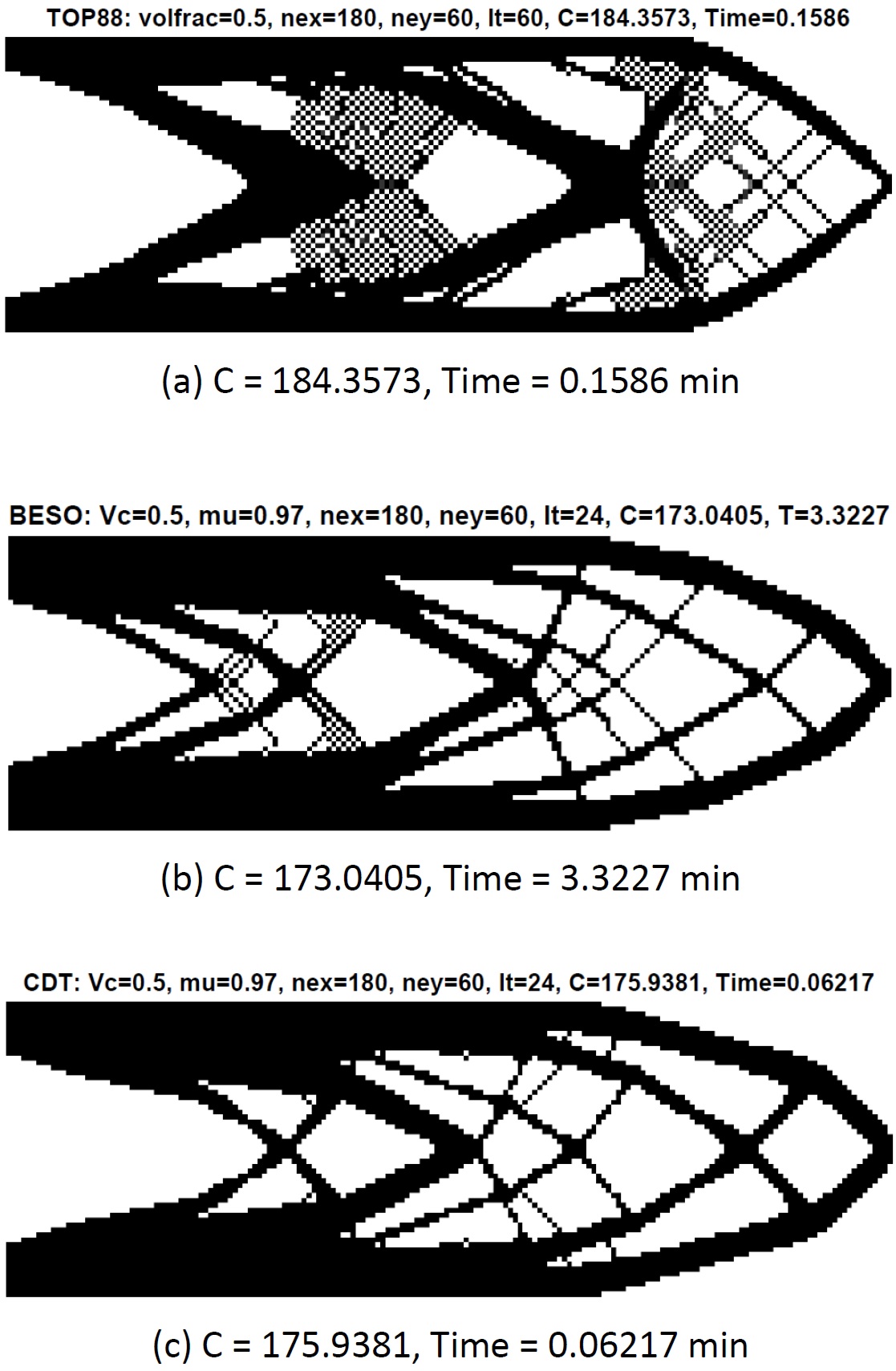}
  \caption{Computational results by SIMP (a), BESO (b) and CDT (c).}
  \label{fig2c}
\end{figure}

 The second example is the 2-D  classical long cantilever beam    (see  Fig.~\ref{fig-cant}).
For this benchmark problem, we first let the volume fraction  $ \mu_c = 0.5$.
Computational  results obtained by the  CDT and by
 SIMP and BESO    are summarized in Fig. \ref{fig2c}.
 Clearly, the  precise solid-void solution  produced by the CDT method    is    much better   than  that by other methods.

\begin{figure} [h]
  \centering
  \includegraphics[width=1\textwidth]{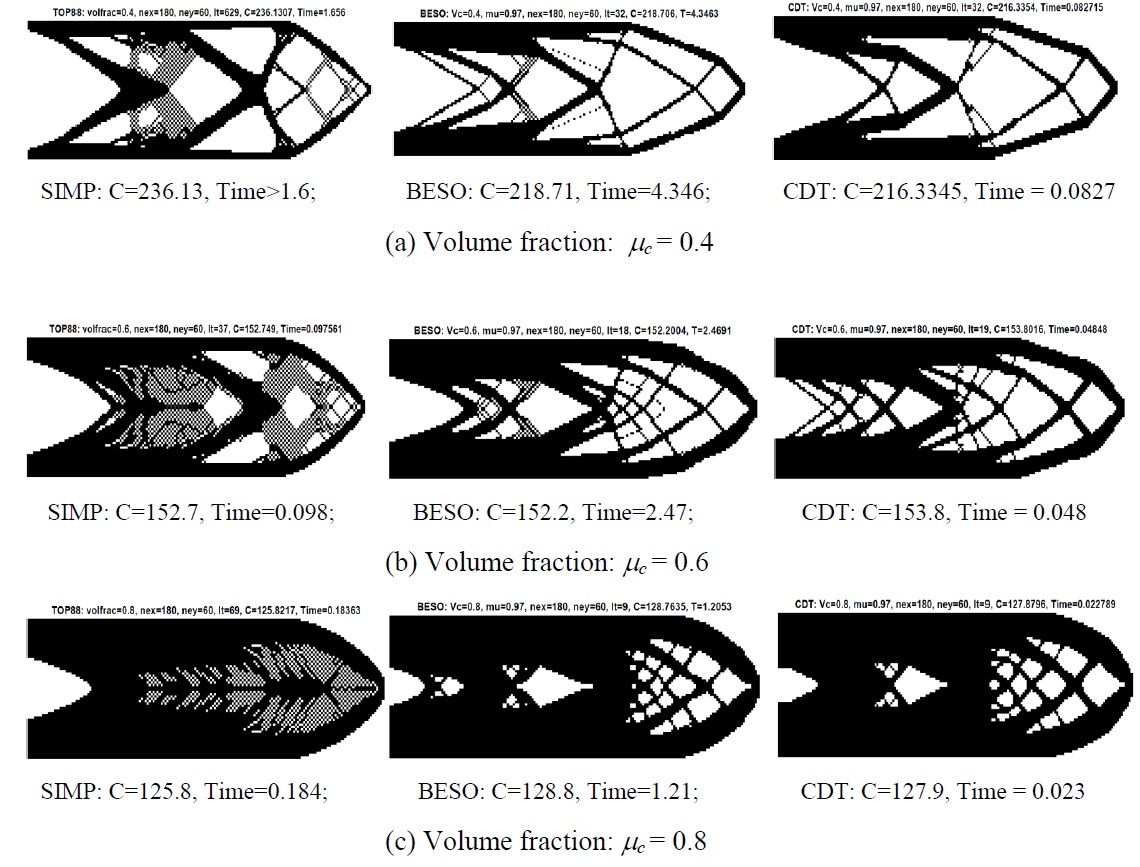}
  \caption{Comparison of computational results by SIMP (left), BESO (middle) and CDT (right) with different volume fractions.}
  \label{fig3c}
\end{figure}

The second test for this problem is the comparison of the three methods with different volume fractions.
 From Fig. \ref{fig3c} one  can see clearly that  the CDT produces   mechanically sound structures with    the shortest computing  time.
A large range of checkerboards is  observed in the results by  SIMP,   a small range of checkerboards is  observed in the results by BESO,  while almost no such pattern  for the proposed CDT method.
For  $\mu_c = 0.4$, the SIMP code does not converge and the result reported in  Fig. \ref{fig3c}(a)  is the output at the 629th iteration.
This test also revealed an important truth, i.e. as the volume fraction $\mu_c$ is decreasing, the strain energies produced by all three methods are  increasing instead of decreasing. This shows clearly that the  minimum strain energy
problem is incorrect for topology optimization.

\begin{figure} [h]
  \centering
  \subfigure[Convergence test for SIMP  (red), BESO (green), and CDT (blue)] {
\includegraphics[width=.6\textwidth, height=48mm]{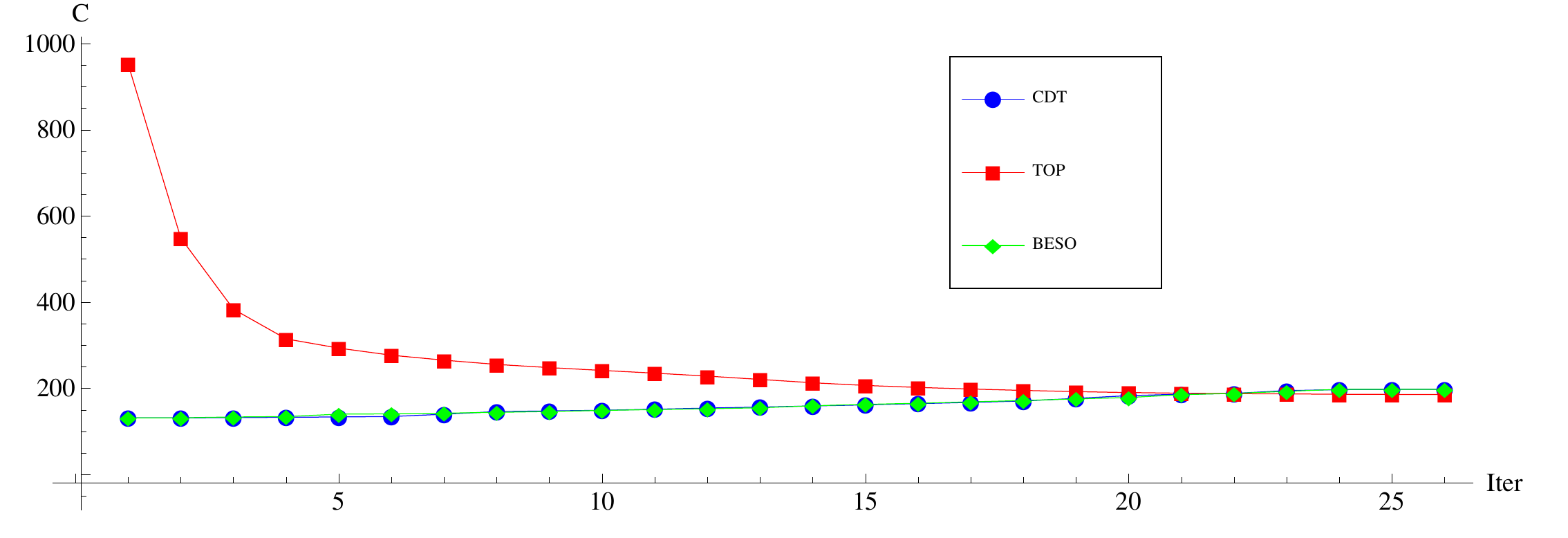}}
  \subfigure[Convergence test for BESO (red) and CDT (blue), $V_c = 0.5, \;\; \mu=0.975$] {
\includegraphics[width=.6\textwidth, height=48mm]{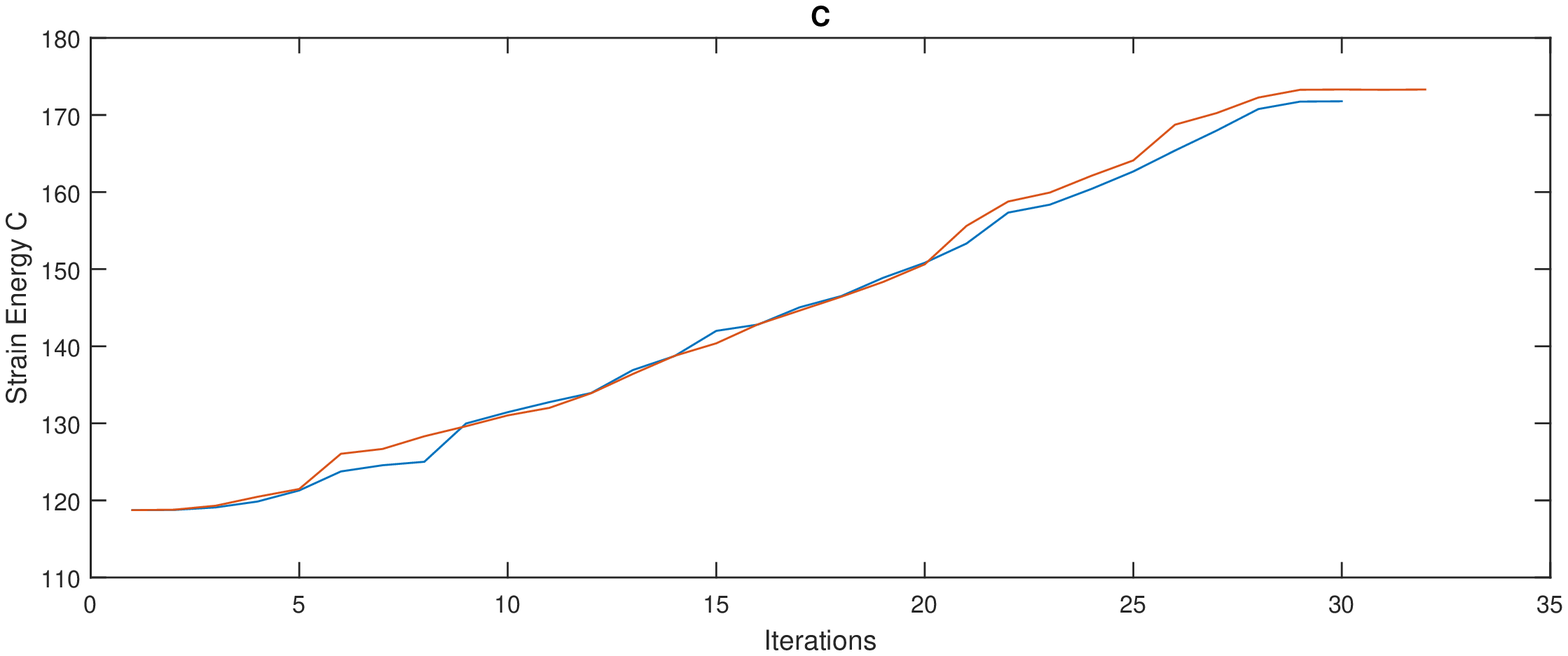}}
\subfigure[Computing times by BESO (red) and CDT (blue), $V_c = 0.5, \;\; \mu=0.975$] {\includegraphics[width=.6\textwidth, height=44mm]{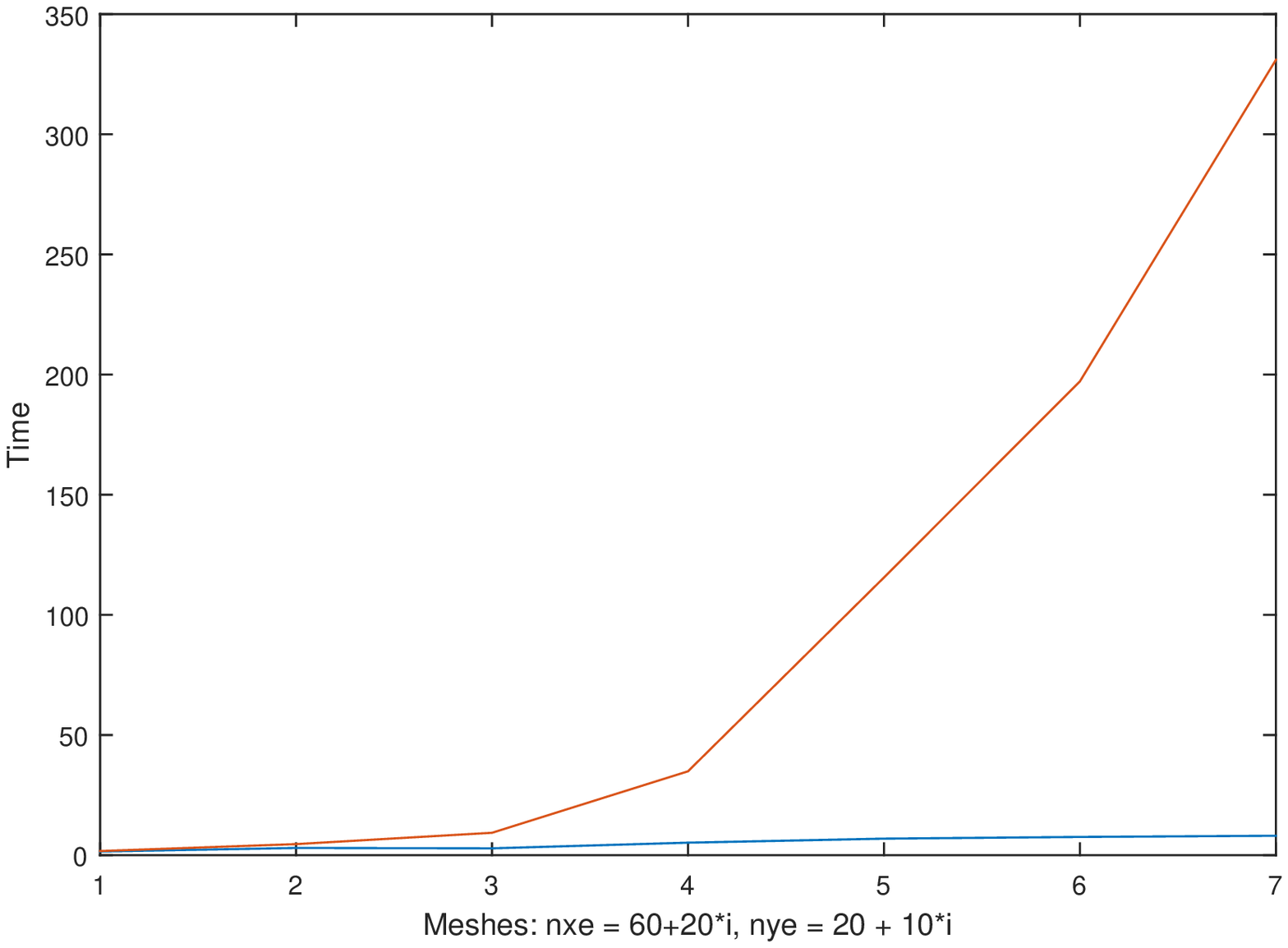}}
 \caption{Comparison tests for SIMP,  BESO and CDT}
  \label{fig-beso-cdt}
\end{figure}

By the fact that the  BESO is a direct method for solving the correct  knapsack problem, although  it can produce the results similar to those by the CDT,  it is not a polynomial-time algorithm.  Fig. \ref{fig-beso-cdt} verified this truth.
For the given mesh resolution $180\times 60$, and $V_c=0.5$, $\mu=0.975$, both methods produce similar target function (i.e. strain energy $C(\brho,\bu) = \half \bu^T \bK(\brho) \bu $)
(see Fig. \ref{fig-beso-cdt}(a,b)), but the BESO's computing time is exponentially  blowing up  as the  increase in  the  mesh numbers (see (see Fig. \ref{fig-beso-cdt}(c)).   Results in Fig. \ref{fig-beso-cdt}(a,b)  show that
 both BESO and CDT are  maximizing   the strain energy, i.e. they are solving the  knapsack problem $(\calP_{kp})$,
while the SIMP method is minimizing  the  strain energy, i.e. it is solving $(\calP_s)$.

\begin{figure} [h]
  \centering
  \includegraphics[width=.6\textwidth]{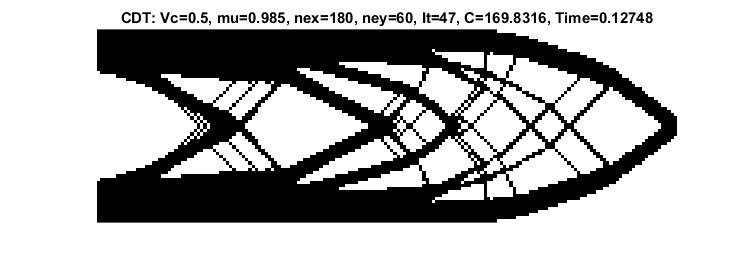}
   \includegraphics[width=.6\textwidth]{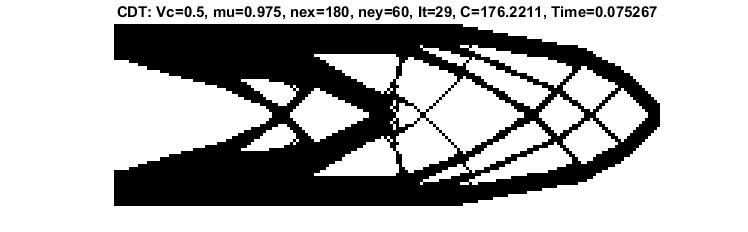}
      \includegraphics[width=.6\textwidth]{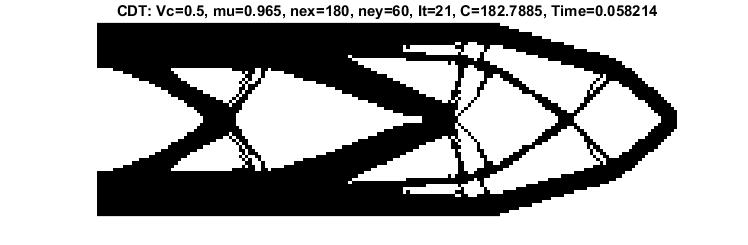}
         \includegraphics[width=.6\textwidth]{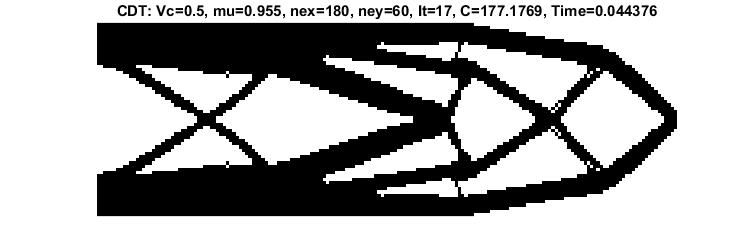}
            \includegraphics[width=.6\textwidth]{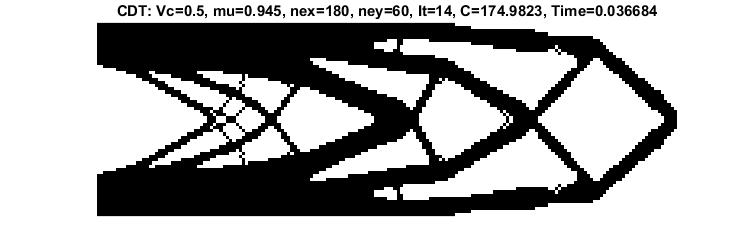}
              \includegraphics[width=.6\textwidth]{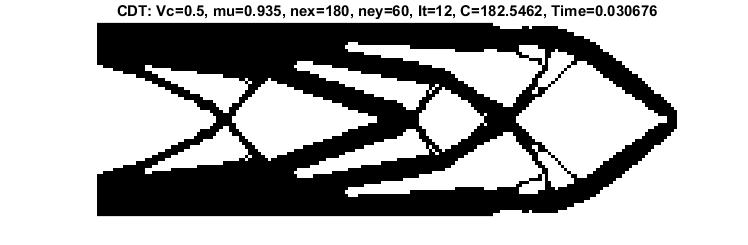}
  \caption{ Experiment for  volume reduction  rate $\mu$ from $ 0.985$ (top) to $ 0.935$ (bottom) }
  \label{fig4c}
\end{figure}

   The reduction  rate $\mu$  in Algorithm~\ref{alg-CDT} plays an important role for both
convergence  and  final result.
We compare in Fig \ref{fig4c} the computational results at different evolutionary rate from  $\mu =0.925$ to $\mu= 0.985$. As we can see
that   big   $\mu$ produces more delicate  structure but  requires more   computational time.
Results in Fig. \ref{fig4c} verify a truth in iteration method  for  bilevel optimization, i.e. the optimal solution depends strongly on the parameter $\mu$.

\begin{figure} [h]
  \centering
  \includegraphics[width=.6\textwidth]{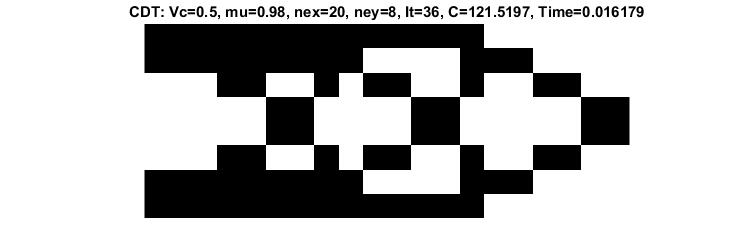}
   \includegraphics[width=.6\textwidth]{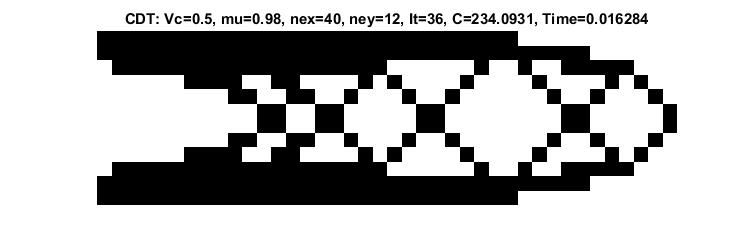}
      \includegraphics[width=.6\textwidth]{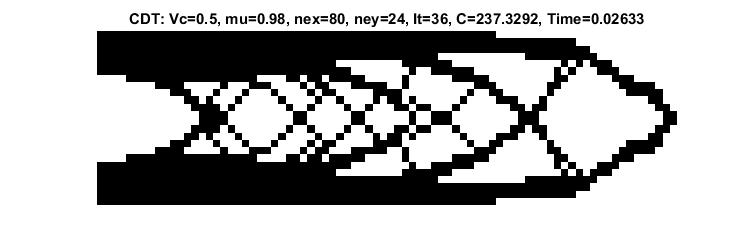}
         \includegraphics[width=.6\textwidth]{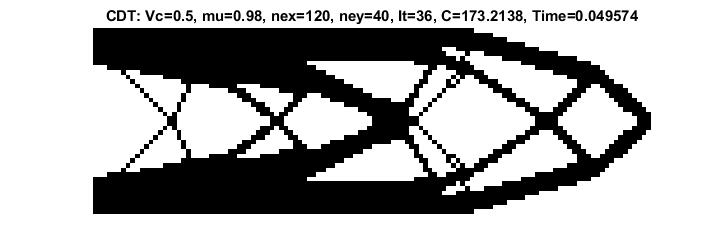}
            \includegraphics[width=.6\textwidth]{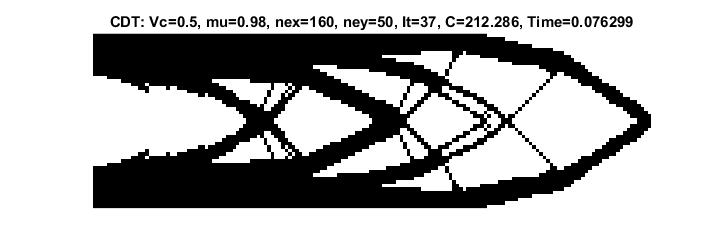}
              \includegraphics[width=.6\textwidth]{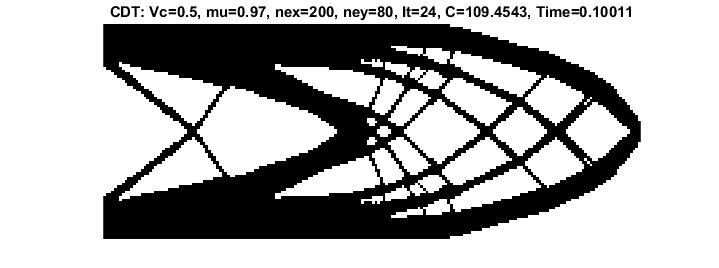}
  \caption{ Experiment for  mesh resolution from $20\times 8$ (top) to $200\times 80$ (bottom) }
  \label{fig4}
\end{figure}

By the fact that the alternative iteration   is adopted in CDT method, the solution also depends sensitively on   mesh resolution.
 Fig \ref{fig4} shows  the CDT solutions to the cantilever beam with different  mesh resolutions from
 $20\times 8$ to $200\times 80$.
 Clearly, different mesh resolution leads to different topology. This is completely reasonable since different mesh size leads to totally different knapsack problem.
 Certainly, the finer is the mesh resolution, the better is the topology.

 \subsection{3-D Cantilever Beam}
 A CDT 3-D code (CDT3D)  is  developed based on the 169 line 3d topology optimization code (TOP3D) by Liu and Tovar \cite{liu-tovar}.
Here we only provide a simple application to the cantilever beam. Since the BESO is not a polynomial-time algorithm, we only
compare the CDT with the SIMP. The parameters in TOP3D are chosen to be penal=3, rmin=1.5.
The volume reduction  rate used in CDT3D is $\mu=0.95$.

\begin{figure}[tb]
  \centering
  \includegraphics[width=0.9\textwidth]{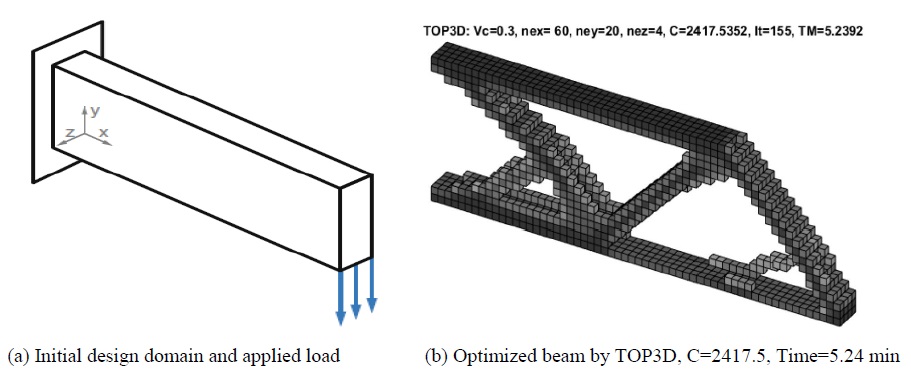}
  \caption{ (a) 3-D cantilever beam and   (b)  SIMP solution }\label{fig6}
 \vspace{0.2in}
\centering
\subfigure[CDT solution]{
\includegraphics[width=1.2\textwidth] {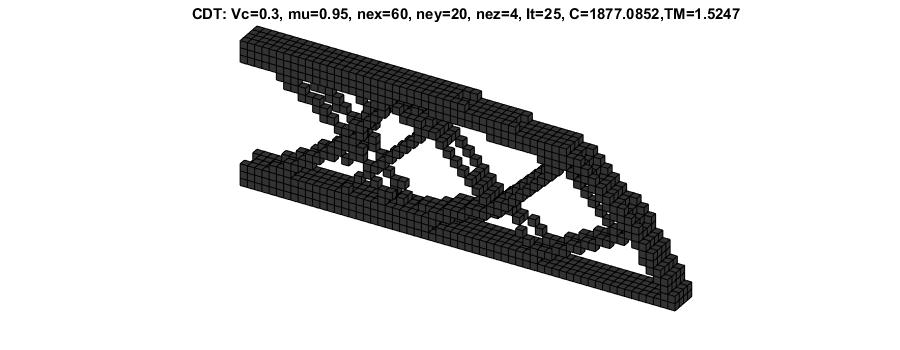} }
\subfigure[Front view ]{
\includegraphics[width=.9\textwidth] {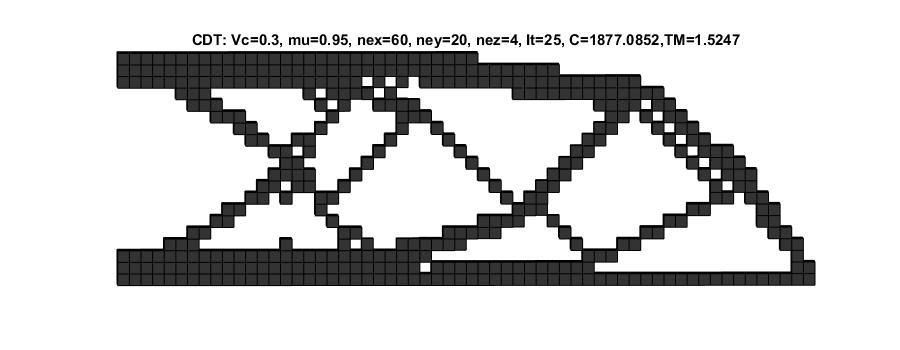} }
\qquad
\caption{CDT solution  and its front view}
\label{fig7}
\end{figure}

For a given volume fraction $\mu_c = 0.3$, Fig. \ref{fig6}(b)  shows the optimized beams by TOP3D.
 Figures  \ref{fig7} shows the solution by CDT3D.

\begin{figure} [h]
  \centering
  \includegraphics[width=.8\textwidth]{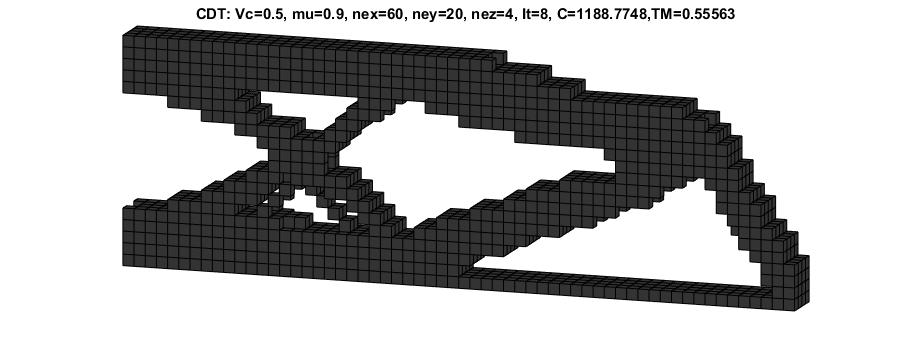}
   \includegraphics[width=.8\textwidth]{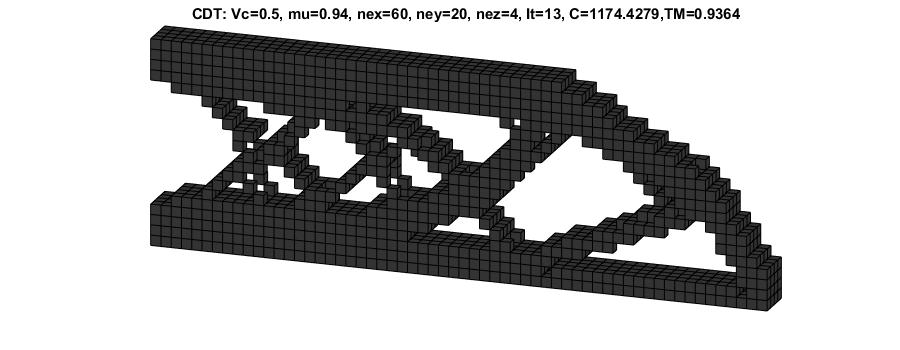}
  \caption{ Experiment for volume reduction  rates  $\mu=0.9$ (top) and  $\mu=0.94$ (bottom) }
  \label{fig8}
\end{figure}

  Fig \ref{fig8} shows  the CDT3D solutions to the cantilever beam with   $\mu=0.9$  and  $\mu=0.94$,
  which
  verified again  a truth in iteration method for solving a coupled optimization problem, i.e. the optimal solution depends strongly on the parameter $\mu$.


\begin{figure}[h]
\centering
\subfigure[ C=21427, Time=12.2]{
\includegraphics[width=.9\textwidth] {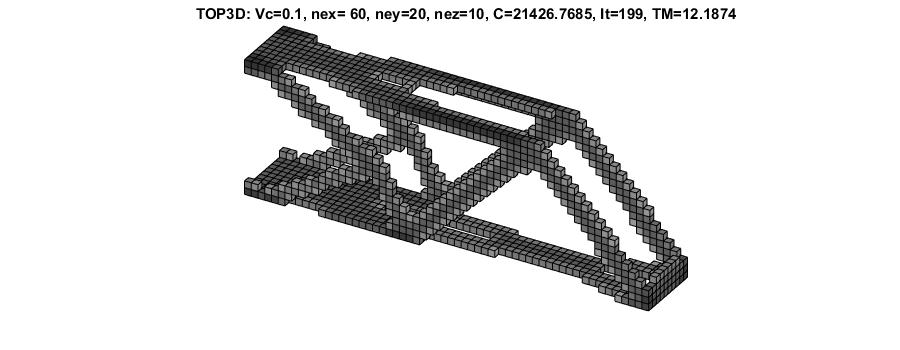} }
\subfigure[Front view of SIMP solution ]{
\includegraphics[width=0.7\textwidth] {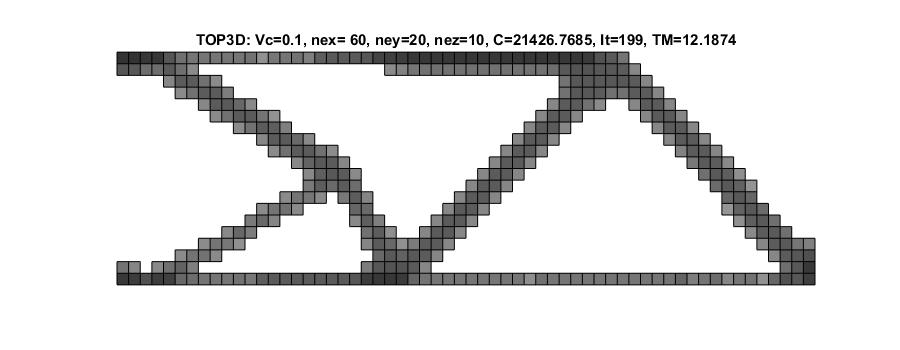} }
\qquad
\caption{SIMP solution with $nez=10$ $\mu_c = 0.1$}
\label{fig9a}
\end{figure}

In order to have a closed look at inside of the 3-D beam, we increase the beam thickness from $nez=4$ to $nez=10$
and decrease the volume fraction from $\mu=0.5$  to $\mu_c = 0.1$.    Fig. \ref{fig9a} and  Fig. \ref{fig9a} show  results   obtained by  the TOP3D  and CDT  (with $\mu=0.93$) methods.
We can see clearly that TOP3D produces many gray elements.  While the optimal structure by the  CDT is elegant and  mechanically sound,
also  the 
 computing time is three times faster.
\begin{figure}[h]
\centering
\centering
\subfigure[C=13193, Time=4]{
\includegraphics[width=.8\textwidth] {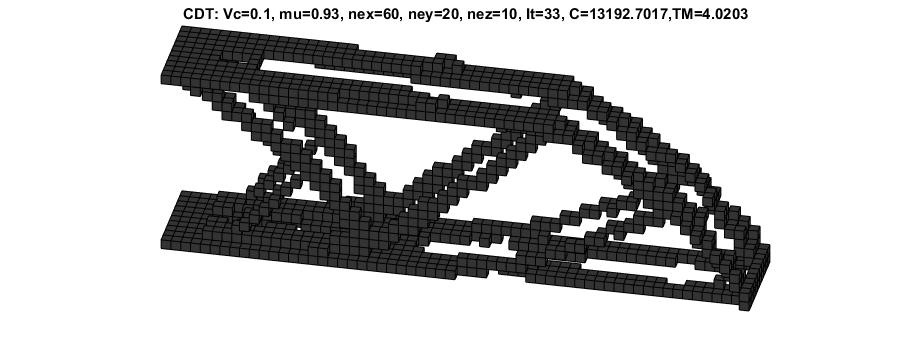} }
\subfigure[Front view of CDT solution ]{
\includegraphics[width=0.7\textwidth] {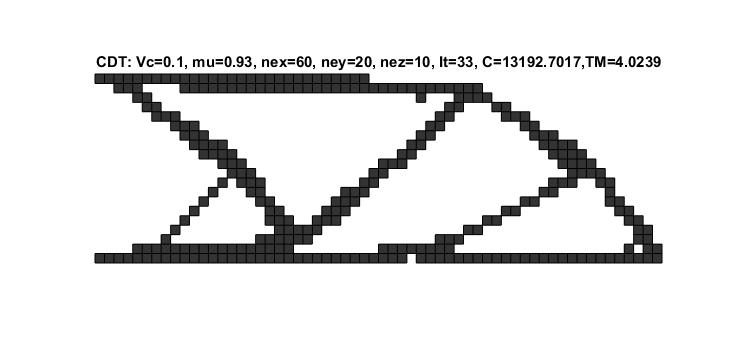} }
\qquad
\caption{SIMP solution with $nez=10$, $\mu_c = 0.1$ and $\mu=0.93$}
\label{fig9b}
\end{figure}

\section{Concluding remarks  and open problems}\label{sec-con}
 Topology optimization has been studied via the canonical duality theory.
Theoretical results presented in this paper   show that this methodological  theory is indeed  powerful
 not only for solving the most challenging topology optimization problems,
but also for correctly understanding and modeling multi-scale problems in complex systems.
The  numerical results  verified  that
the CPD  method can produce mechanically sound optimal topology, also it is   much more faster than the popular SIMP and BESO methods.
Specific conclusions are given  below.
\begin{enumerate}
 \item The mathematical problem  for general topology optimization should be formulated as a  bi-level
knapsack  programming $(\calP_{bk})$. This model works for both linearly and nonlinearly deformed elasto-plastic structures.

 \item   The CDT  is a polynomial-time algorithm, which can solve  general topology optimization problem to obtain global optimal solution at each volume  iteration.

\item The  $\beta$-perturbation  for solving integer programming problems is actually the canonical penalty-duality method  proposed in \cite{gao-cs88}.
Theorem \ref{thm-solution}  is  an application of the pure complementary energy principle in nonlinear elasticity \cite{gao-mrc99}.
 This principle plays an important role not only in nonconvex analysis and  computational mechanics, but also in  topology optimization, especially for large deformed structures.

 \item  Unless a magic method is proposed, the volume reduction method  is necessary for solving  general knapsack-type problems
  if $ \mu_c = V_c/V_0 \ll 1$. But the global optimal solution  depends  sensitively on the reduction  rate $\mu \in [\mu_c, 1)$.

\item The   so-called   compliance  minimization problem $(\calP_c)$   is actually a single-level reduction for solving $(\calP_{bk})$.
 Alternative iteration for  the  strain energy  minimization $(\calP_s)$ leads to an anti-knapsack problem, which has only trivial solution.

\item The   problem $(\calP_{simp})$  is a box-constrained nonconvex minimization subjected to a knapsack condition.
The  SIMP    is not a mathematically correct penalty method for solving either $(\calP_c)$ or  $(\calP_s)$.
Even if the  magic number $p=3$ works  for certain materials,  this  method can't produce  correct
integer  solutions for any given $p > 1$.

\item The  BESO  method  is actually  a combination of the volume reduction and a direct method for solving the  knapsack problem $(\calP_{kp})$.
For small-scale problems,  the  BESO    method  can produce  reasonable results  better  than by SIMP.
But  it is  not a  polynomial-time algorithm.
\end{enumerate}

 By the fact that the general topology optimization  problem $(\calP_{bk})$ is a coupled, mixed integer nonlinear programming,
  although  for each given lower-level solution $\bu\in \calU_a$  the  upper-level problem  can be solved analytically by the canonical duality theory,
  the final result  produced by the CDT may not be the global optimal solution to $(\calP_{bk})$
  since the alternative iterations for $\brho$ and $\bu$ are adopted
 in Algorithm 1.
This is the main reason why the numerical solution produced by the CDT depends on mesh resolution and volume reduction   rate $\mu$.
The main open  problems include     the optimal  parameter
  $\mu$  in order to ensure the fast convergence rate with  the optimal results,  the influence of the mesh resolution to the optimal topology,  the existence and uniqueness of a  global optimization solution
  to $(\calP_{bk})$ for a given design domain $V_c$.

This paper is just  a simple  application of the canonical duality theory for  linear elastic topology optimization.
A  66-line Matlab code for topology optimization will be  available soon to public. Also a refined algorithm with a 44-line Matlab code
has been developed.
 The  canonical duality  theory is particularly useful for studying  nonconvex, nonsmooth, nonconservative  large deformed dynamical systems \cite{gao-royal}.
 Therefore, the future works include
 the CDT method for solving  general topology optimization problems of  large deformed elasto-plastic  structures subjected to dynamical loads.
\vspace{-.5cm}

\section*{Appendix: A Brief Review on Canonical Duality Theory }
\vspace{-.5cm} The canonical duality theory
 proposed in \cite{gao-dual00} comprises mainly
 \begin{verse}
  i) a {\em canonical dual transformation}, \\
   ii) a
{\em complementary-dual principle,} and \\
 iii) a {\em triality theory. }
\end{verse}
 The canonical dual transformation is a versatile
methodological method  which can be used to model complex systems within a unified framework
and to formulate perfect dual problems without a duality gap.
The complementary-dual principle presents a unified analytic solution form for general problems in continuous and
discrete systems. The triality theory reveals an intrinsic duality pattern in multi-scale  systems, which can be used
to identify both global and local extrema, and to develop
deterministic  algorithms for effectively solving a wide class of
nonconvex/nonsmooth/discrete optimization/variational problems.

The canonical duality theory was developed from Gao and Strang's original work \cite{gs-89}  for solving the
following   nonconvex/nonsmooth variational  problem
\eb
\inf \left\{ \Pi(\bu) = \int_\Omega W(\nabla \bu) \dO - F (\bu) \;  | \;\; \bu \in \calU   \right\}  ,  \label{eq-gv}
\ee
where $F(\bu)$ is a generalized  external  energy, which is  a linear functional of $\bu$
 on its effective domain $\calU_a \subset \calU$;
  the internal energy density
 $W(\bF)$ is  the so-called {\em super-potential} in the sense of J.J. Moreau \cite{moreau}, which is a nonconvex/nonsmooth function of the deformation gradient $\bF = \nabla \bu$.
Numerical discretization of  (\ref{eq-gv})   leads to a global optimization problem in finite dimensional space  $\calX \subset \real^n$.
 Canonical dual finite element method   for solving the nonconvex/nonsmooth mechanics problem (\ref{eq-gv}) was first proposed   in 1996 \cite{gao-jem96}.

\begin{remark}[Objectivity and Conceptual Mistakes]  {\em
Objectivity is a central  concept in our daily life, related to reality and truth. In philosophy, it   means the state or quality of being true even outside a subject's individual biases, interpretations, feelings, and imaginings\footnote{\url{https://en.wikipedia.org/wiki/Objectivity_(philosophy)}}.
In science, the objectivity  is often attributed to the property of scientific measurement,
as the accuracy of a measurement can be tested independent from the individual scientist who first reports it\footnote{ \url{https://en.wikipedia.org/wiki/Objectivity_(science)}}.
In continuum mechanics, a real-valued function $W :\calF \subset \real^{n\times m}  \rightarrow \real $ is called to be {\em objective}  if
\eb
W(\bF) = W(\bR \bF) \;\; \forall  \bF  \in \calF, \;\; \forall \bR, \;\;  \bR^T = \bR^{-1}, \;\; \det \bR = 1.
\ee
The objectivity  plays a fundamental rule in mathematical modeling, which is also called the principle of material frame indifference \cite{oden}.

However,  this
 important concept    has been misused in
    optimization literature. As a  result,  the general problem in  the so-called nonlinear  programming  is  formulated as
\eb
\min f(x), \;\; s.t. \;\; h(x) = 0, \;\;  g(x) \le 0, \label{eq-ap}
\ee
where $f(x)$ is called the ``objective function'',\footnote{This terminology is used mainly  in English literature.
It  is correctly called  the target function in all Chinese and Japanese literature,  or the goal function in some   Russian and German literature.}  which
is allowed to be any arbitrarily given function,  even a  linear function.
Clearly, this mathematical problem is too abstract.
Although it enables one to  model  a very wide range of problems, it comes at a price: many
global optimization problems are considered to be NP-hard.
Without detailed information on these arbitrarily given functions, it is impossible to have a powerful theory for solving the
artificially given constrained  problem \eqref{eq-ap}. \hfill $\clubsuit$ }
  \end{remark}

In linguistics, a complete and grammatically correct sentence should be composed by at least three words:  subject, object, and a predicate.
As a language of science,  the mathematics should follow this rule. Based on the excellent works by
Oden-Reddy, Strang\footnote{The celebrated textbook \cite{strang}  by Gil Strang is a required course for all engineering graduate students
 at MIT. Also, the well-known  MIT online teaching program was started from  this course.}, and Tonti \cite{oden,oden-reddy,strang,tonti},
 a general mathematical problem  was proposed in  \cite{gao-aip}:
  \eb
 (\calP_g): \;\;   \min \{\Pi (\bx) =  W (\bD  \bx) - \FF (\bx) \;  | \;\; \bx \in \calX \} , \label{eq-go}
   \ee
   where $\bD $ is a linear operator ($\nabla$ is a special case),  which maps  the configuration space  $ \calX$ into a so-called intermediate space ${\cal W}$ \cite{oden,oden-reddy,tonti};
   $W:{\cal W} \rightarrow \real \cup \{ + \infty\}$ is an objective function such that its sub-differential leads to
   the {\em  constitutive duality}  $ \bw^* \in \partial^- W(\bw) \subset { \cal W}^*$, which is  governed by
   the constitutive law and all possible physical constraints;
 while   $ \FF:\calX \rightarrow \real\cup\{ -\infty\}$ is a so-called {\em subjective function} such that its over-differential
   leads to the {\em action-reaction duality}  $\bx^* \in \partial^+  \FF (\bx) \subset \calX^*$  including all possible boundary conditions and geometrical constraints \cite{gao-jogo00}.
   This subjective function must be linear on its
   effective domain  $\calX_a  = \{ \bx \in \calX| \;\; \FF(\bx) > - \infty \}$, i.e. $ \FF (\bx) = \la  \bx ,  \bff \ra \;\; \forall \bx \in \calX_a$, where  $\bff \in \calX^*$ is an given action,
   the feasible set $\calX_a$    includes all possible geometrical constraints.
   The predicate in $(\calP_g)$ is the operator ``$-$'' and the difference $\Pi(\bx)$ is  called the target function in general problems.
The object and subject are in balance only at the optimal states.
  By the fact that $\calX$ and $\calW$ can be in different dimensional spaces, this unified problem can be used to model multi-scale systems \cite{gao-dual00,gao-yu}.

  The unified  form  $(\calP_g)$  covers general  constrained nonconvex/nonsmooth/discrete  variational and optimization problems in  multi-scale complex systems
   \cite{g-l-r-17,gao-yu} as well as all equilibrium problems in  linear elasto-plasticity and mathematical physics \cite{oden-reddy,strang}.

  \begin{definition}[Properly and Well-Posed Problems]
  A  problem   is called {\em properly posed} if
for any given non-trivial input it has  at least one non-trivial solution.
It is called {\em well-posed} if the solution is  unique.
\end{definition}
Clearly, this definition is more general than Hadamard's well-posed problems in dynamic  systems
since   the continuity  condition is not included.
 Physically speaking, any real-world problems  should be well-posed since
  all natural phenomena  exist   uniquely.
  But practically, it is difficult to model  a real-world problem precisely.
  Therefore,  properly posed problems are allowed for the canonical duality theory.
  This definition is important for understanding the triality theory and NP-hard problems.

 It is well-known in finite deformation theory that the stored energy $W(\bF)$ is an objective function if and only if
  there exists  an objective strain tensor
 $\bC = \bF^T \bF$  and a  function $U(\bC)$  such that $W(\bF ) = U(\bC)$ \cite{ciarlet}.
  In nonlinear elasticity, the function $U(\bC)$  is usually convex (say the St Venant-Kirchhoff material) such that  the duality relation
  $\bC^* \in  \partial U(\bC)$  is invertible and the
 complementary energy $U^*(\bC^*)$ can be uniquely defined via the
 Legendre-Fenchel transformation $U^*(\bC^*) = \sup_{\bC}  \{ \bC: \bC^* - U(\bC) \}$.
   These basic truths in continuum physics laid a foundation for the canonical duality theory.

   The key idea of  the  {\em canonical transformation}    is to  introduce a nonlinear  operator
   $\bxi = \Lam (\bx) :\calX \rightarrow \calE$ and a convex, lower semi-continuous function $\Psi(\bxi): \calE \rightarrow \real \cup \{ + \infty\}$ such that
  the nonconvex function $W (\bD \bx)$ can be written in the canonical form
  \eb
  W (\bD \bx) = \Psi(\Lam(\bx)) \;\; \forall \bx \in \calX
  \ee
  and the following canonical duality relations  hold
   \[
   \bzeta \in \partial \Psi(\bxi) \;\; \Leftrightarrow \;\; \bxi \in \partial \Psi^*(\bzeta) \;\;  \Leftrightarrow \;\;
   \Psi(\bxi) + \Psi^*(\bzeta) = \langle  \bxi ; \bzeta \rangle ,
   \]
 where $\langle  \bxi ; \bzeta \rangle $ represents the bilinear form on $\calE$ and its canonical dual $\calE^*$.
   Thus, using the Fenchel-Young equality $   \Psi(\bxi)  = \langle  \bxi ; \bzeta \rangle -  \Psi^*(\bzeta)$,
   the nonconvex/nonsmooth minimization problem (\ref{eq-go}) can be equivalently written in the following min-max form
   \eb
   \min_{\bx \in \calX} \max_{\bzeta \in \calE^*} \{ \Xi(\bx, \bzeta) = \langle \Lam(\bx) ; \bzeta \rangle - \Psi^*(\bzeta) -  \FF ( \bx)  \},
   \ee
where $ \Xi(\bx, \bzeta)$ is the so-called  {\em total complementary function}, first
introduced by Gao and Strang in 1989 \cite{gs-89}.
Particularly, if  $\bveps= \Lam(\bx)$ is a homogeneous quadratic  operator, then
$  \langle \Lam(\bx) ; \bzeta \rangle =  \half  \bx^T \bG(\bvsig) \bx  = G_{ap}(\bx, \bvsig)$ is the so-called {\em complementary
  gap function}, where   the Hessian matrix $\bG(\bvsig) = \nabla^2_{\bx} G_{ap}(\bx, \bvsig)$ depends on the quadratic operator $\Lam(\bx)$.
  This gap function  recovers
the duality gap  produced by the traditional Lagrange multiplier method and  modern Fenchel-Moreau duality theory.
By using $\Xi$, the canonical dual  function  $\Pi^d(\bvsig)$ can be obtained  as \cite{gao-cace}
\eb
\Pi^d  (\bvsig)  = \mbox{sta}  \{ \Xi(\bx, \bvsig) | \;\; \bx \in \calX_a \;  \} = - \half  \la    [\bG(\bvsig)]^{+} \bff , \bff \ra    - \Psi^*(\bvsig),
\ee
where $\mbox{sta } \{ f(x) |\; x \in X \}$ stands for finding a stationary value of $f(x),\ \forall x\in X$,
and  $\bG^+$ represents a generalized inverse of $\bG$.
Based on this canonical dual function, a  general  solution to the nonconvex   problem  (\ref{eq-go})  can be obtained \cite{gao-jogo00}.
\begin{thm}[Complementary-Dual Principle and Analytic Solution]
$\;$  For a given input  $\bff \in \calX^*$, the vector
 $\barbx = [\bG(\barbvsig)]^{+} \bff $
is a  stationary point of $\Pi  (\bx)$ if and only if
  $\barbvsig$ is a stationary point of $\Pi^d  (\bvsig)$,   $(\barbx, \barbvsig)$ is a stationary point
   of $\Xi(\bx, \bvsig)$,  and
$ \Pi  (\barbx) = \Xi(\barbx, \barbvsig) = \Pi^d  (\barbvsig) $
\end{thm}
This theorem was first proposed in
     1997 for a large deformed beam theory \cite{gao-amr97}.
     In nonlinear elasticity, the canonical dual function $\Pi^d  (\bvsig) $ is the so-called pure complementary energy, first proposed in 1999 \cite{gao-mrc99}.
This  principle   solved a 50-year old open problem in finite deformation theory
 and is known as the Gao principle in literature \cite{li-gupta}.
 Clearly, $\bvsig \in \calS^+_a$ if and only if $G_{ap}(\bx, \bvsig) \ge 0 \;\; \forall \bx\in \calX_a$.
The following triality  theory shows that this  gap function  provides a global optimality criterion in nonconvex mechanics and global optimization.
\begin{thm}[Triality Theory \cite{gao-jogo00}]\label{thm-tri}
Suppose that    $(\barbx, \barbvsig)$ is a stationary solution of $\Xi(\bx, \bvsig)$ and   $\calX_{o}\times \calS_{o} $  is  a neighborhood of $\left(
\barbx,\barbvsig \right) $.

If  $\barbvsig \in \calS^+_a = \{ \bvsig \in \calE^* | \;\; \bG(\bvsig )  \succeq 0 \}$, then
\vspace{-.1cm}
\begin{equation}
\Pi  (\barbx) = \min_{\bx \in \calX_a}  \Pi   \left( \bx\right)  =\max_{\bvsig \in \calS^+_a}  \Pi^{d}  \left(\bvsig \right) = \Pi^d  (\barbvsig) . \label{Global0} \vspace{-.1cm}
\end{equation}

If $\barbvsig \in \calS^-_a = \{ \bvsig \in \calE^* | \;\; \bG(\bvsig )  \prec 0 \}$,
then  we have either  \vspace{-.2cm}  
\begin{equation}
\Pi (\barbx) = \max_{\bx \in \calX_o}  \Pi  \left( \bx\right)  =\max_{\bvsig \in \calS_o}  \Pi^{d}  \left(\bvsig \right) = \Pi^d  (\barbvsig)  ,  \label{dmax} \vspace{-.1cm}
\end{equation}
or (if $\dim \Pi  = \dim \Pi^d  $) \vspace{-.1cm}  
\begin{equation}
\Pi  (\barbx) = \min_{\bx \in \calX_o}  \Pi   \left( \bx\right)  =\min_{\bvsig \in \calS_o}  \Pi^{d}  \left(\bvsig \right) = \Pi^d (\barbvsig) .  \label{dmin} \vspace{-.1cm}
\end{equation}
\end{thm}

\begin{figure}[h]
  \centering
  \subfigure[$f=0.5$]{\includegraphics[width=0.48\textwidth]{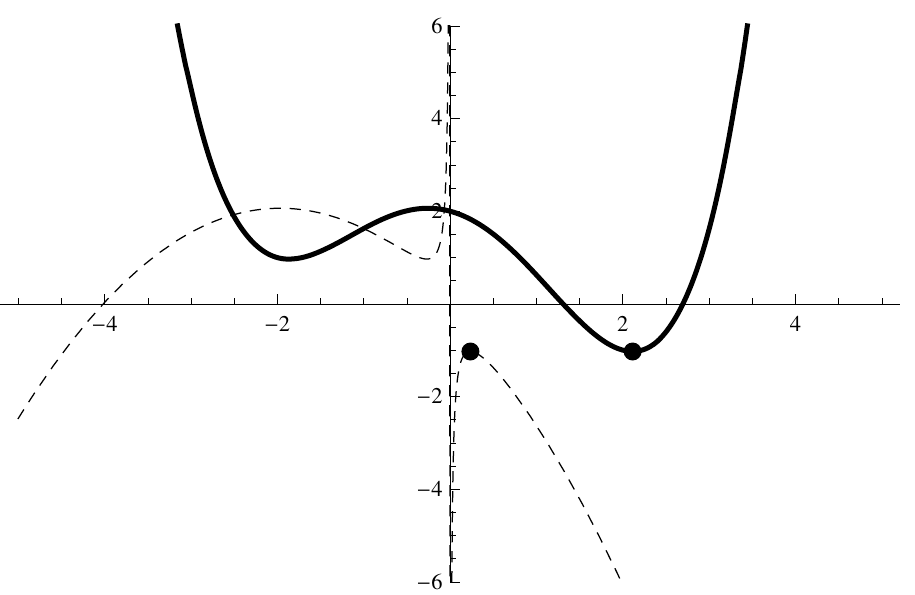}}\quad
  \subfigure[$f=0$]{\includegraphics[width=0.48\textwidth]{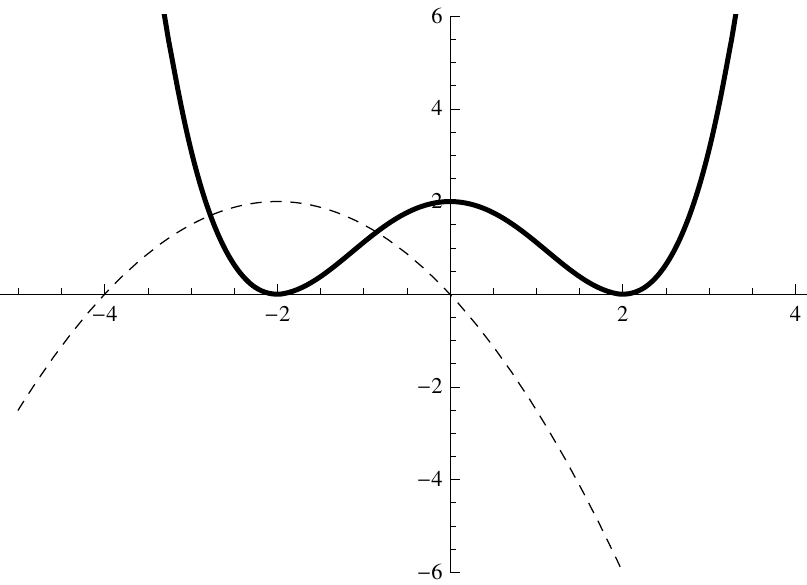}}
\caption{ Graphs of   $ \Pi (\bx)$  (solid) and $\Pi^d(\vsig)$ (dashed) for $n=1, \;\; \beta = 1, \; \; \lam= 2$ with $x_1= 2.1$, $\vsig_1=0.24$ and $\Pi(x_1)=-1.02951=\Pi^d(\vsig_1)$ (see the two black dots in (a)) if $f=0.5$. While $x_{1,2} = \pm 2$ if $f=0$ (see (b)).     }
\label{fig-dw}
\end{figure}

This triality theory shows that the Gao-Strang gap function can be used to identify both global minimizer and the biggest local extrema.
To understand the canonical duality theory, let us consider a simple nonconvex optimization in
$\real^n$:
\eb
\min \left\{  \Pi(\bx)=\half \beta (\half\|\bx\|^2-\lam)^2-\bx^T \bff  \;\;\;   \forall \bx \in \real^n \right\} ,
\ee
where $\alp, \lam > 0$ are given parameters.
The criticality condition $\nabla \Pi (\bx)=0$ leads to a nonlinear algebraic
equation system  $
 \beta  (\half \|\bx\|^2-\lam)\bx =\bff.
$
Clearly,  this n-dimensional nonlinear algebraic equation may have multiple solutions.
Due to the nonconvexity,  traditional convex optimization theories
 can't be used to identify global minimizer.  However, by the
canonical dual transformation, this problem can be solved completely to obtain all extremum solutions.
To do so, we let  $\xi=\Lam(\bx)=\half\|\bx\|^2  \in \real$. Then,
the nonconvex function $W(\bx) = \half  \beta  (\half \| \bx \|^2 -\lam)^2$
can be written in canonical form
$\Psi(\xi) = \half \beta (\xi - \lam)^2$.
Its Legendre conjugate is  given by
$\Psi^*(\vsig)=\half  \beta^{-1}\vsig^2 + \lam \vsig $, which is strictly convex.
Thus,
the total  complementary function for this nonconvex optimization
problem is
\eb
\Xi(\bx,\vsig)= \half \|\bx\|^2 \vsig  - \lam\vsig-\half
 \beta^{-1}\vsig^2 - \bx^T \bff.
\ee
For a fixed $\vsig \in \real$, the  criticality condition
$\nabla_{\bx} \Xi(\bx,\vsig)=0$ leads to
$\vsig \bx-\bff=0. $  So the analytical solution
  $\bx= \vsig^{-1} \bff $ holds for $\vsig\neq 0$. Substituting this into the
total complementary function $\Xi$,
the canonical dual function can be easily obtained as
\eb
\Pi^d(\vsig) = \{\Xi(\bx,\vsig)| \nabla_{\bx} \Xi(\bx,\vsig)
=0\} =  -\frac{ \|\bff \|^2  }{2 \vsig}-\half \beta^{-1} \vsig^2
-\lam \vsig, \;\;\; \forall\vsig\neq 0.
\ee
The critical point of this canonical function is obtained
by solving the following dual algebraic  equation
\eb
( \beta^{-1} \vsig+\lam)\vsig^2=\half \|\bff \|^2 . \label{eq-deuler}
\ee
For any given parameters $\beta$, $\lam$ and the vector $\bff  \neq {\bf 0} \in \real^n$,
this cubic algebraic equation has  three non trivial real solutions
  $\vsig_ 1 \in \calS^+_a = \{\vsig \in \real\; |\; \vsig <   0 \}$,  $\vsig_{2,3} \in \calS^-_a = \{ \vsig \in \real\; |\; \vsig < 0 \}$,
then triality theory  tells us that $\bx_1 = \bff/\vsig_1$ is
a global minimizer, $\bx_2  = \bff/\vsig_2$ is a local minimizer, $\bx_3  = \bff/\vsig_3$ is a local maximizer  of $\Pi(\bx)$, and $\Pi(\bx_i) = \Pi^d(\vsig_i), \; i=1,2,3$  (see Fig \ref{fig-dw}).

If  $\bff= 0$, the graph of $\Pi(\bx)$ is symmetric (i.e. the so-called double-well potential or the Mexican hat for $n=2$
\cite{gao-jogo00})
with  infinite number of global minimizers satisfying $\| \bx \|^2  = 2 \lam$.
In this case, the canonical dual $\Pi^d (\vsig) = - \half \beta^{-1}\vsig^2 - \lam \vsig$ is strictly concave
with only one critical point (local maximizer) $\vsig_3 = -  \beta \lam  < 0 $ (for $\beta, \lam > 0$).
The corresponding solution $\bx_3 = \bff \vsig_3 = 0$ is a local maximizer and $\Pi(\bx_3) =\Pi^d(\vsig_3)$.
 By the canonical dual equation (\ref{eq-deuler})
we have $\vsig_1 =\vsig_2 = 0$ located on the boundary of $\calS^+_a$. In this case, the analytical solution  $\bx = \bff /\vsig $ is singular.
For  $n=1$, these two roots are corresponding to the two global minimizers   $x_{1,2} = \pm \sqrt{2 \lam}$ (see Fig. \ref{fig-dw} (b)), which can be obtained by linear perturbation.

This simple example reveals  a fundament truth  in global optimization, i.e.,  for a  given nonconvex/discrete minimization problem
 in $ \calX  \subset \real^n$, its canonical dual is a
concave maximization over a convex set  $ \calS^+_a  \subset \real^m$,
which can be solved easily to obtain a unique solution if $\calS^+_a$ has no-empty interior. Very often  $n \gg m$.
Otherwise, the primal problem could be really NP-hard. Therefore, the canonical duality theory can be used to identify if a global optimization problem is not NP-hard.
Generally speaking, if a non-symmetrically   distributed external force  $\bff \neq 0 $, the canonical dual feasible set $\calS^+_a$ has no-empty interior.
It was proved in \cite{gao-cace} that for quadratic 0-1 integer programming problem, if the
source term $\bff$ is bigger enough, the solution is simply $\{ x_i\} = \{0  \mbox{ if }  f_i  < 0 , \; 1 \; \mbox{ if } f_i > 0\}$ (Theorem 8, \cite{gao-cace}).
 The canonical duality theory has been used successfully for solving a large class of challenging problems in multidisciplinary fields of engineering and sciences.
 A comprehensive review and applications are given in the recent book \cite{g-l-r-17}.

 Mathematics and mechanics have been complementary partners since the Newton times.
Many fundamental ideas, concepts, and mathematical methods extensively used in calculus
of variations and optimization are originally developed from mechanics.
The canonical duality theory is a typical example  developed from Gao and Strang's work on finite  deformation theory, where the subjective function $F(\bu)$ must be linear so that the input force
$\bff = \partial F(\bu)$ is independent of $\bu$ (i.e. a dead-load).
Dually, the objective function $W(\beps)$ must be nonlinear such that there exists an objective measure $\bxi = \Lam(\bu)$ and a  convex function $\Psi(\bxi)$, the
canonical transformation $W(\bD \bu) = \Psi(\Lam(\bu))$ holds for most real-world  systems.
This is  the reason why the canonical duality theory  can be used naturally to solve  general  challenging problems in  multidisciplinary fields.
 However,   the  objectivity has been misused   in engineering  optimization. As a consequence,  the subjective function   $F(\bu) = \bu^T \bff$
  has been extensively   used as the objective function in topology optimization,   the minimum compliance problem $(\calP_c)$ has been  mistakenly written as the minimum strain energy problem $(\calP_s)$,  and the stored energy
   $ W(\bD\bu) = \half  \bu^T \bK (\brho)\bu$  has been confusedly   called  the mean compliance.
Also,  the objectivity  has   been misused   in mathematical optimization. It turns out  the canonical duality    theory has been mistakenly challenged
 by M.D. Voisei and  C. Z$\breve{a}$linescu
in several publications    (cf.  \cite{g-l-r-17}). By oppositely choosing  linear functions for $W(\beps)$ and nonlinear functions for $F(\bu)$,
they produced a list  of ``count-examples''  and concluded: ``a correction of this theory is impossible without falling into trivial''.
 The  conceptual  mistakes in  these   challenges as well as in topology optimization
revealed at least two important  issues: 1) there exists a big  gap between optimization and mechanics;
2)  incorrectly using the well-defined concepts can lead to not only  ridiculous arguments, but also wrong
mathematical model.
As V.I. Arnold indicated \cite{anorld}: ``In the middle of the twentieth century it was attempted to divide physics and mathematics. The consequences turned out to be catastrophic''.
 Interested readers are
recommended to read the recent papers    \cite{gao-opl16,g-r-l}
  for further discussion.
\section*{Acknowledgements}
The author would like to sincerely acknowledge the important comments and suggestions from an anonymous reviewer. A new section  \ref{sec4}  is added in
 the  revised  version to address some fundamental issues in computational complexity, which  is   important for correctly understanding the NP-hardness and challenges   in topology optimization, global optimization and computer science.
 This research is supported by US Air Force Office for Scientific Research (AFOSR)  under the grants   FA2386-16-1-4082 and FA9550-17-1-0151.
The Matlab code was helped by Professor M. Li from Zhejiang University and  Ms Elaf Ali at Federation University Australia.
Both the algorithm and code have been  tested  thoroughly by Dr. Flex Ospald from Chemnitz University of Technology.

\bibliographystyle{plain}
\bibliography{cdtor1}
\end{document}